\documentclass[12pt]{amsart}
\usepackage{amscd,amsthm,amssymb,amsfonts,amsmath,epsfig,epic,eclbip,bbm,url}
\usepackage[arrow,matrix,tips,2cell]{xy}

\theoremstyle{plain}
\newtheorem{theorem}{Theorem}[section]
\newtheorem{lemma}[theorem]{Lemma}
\newtheorem{proposition}[theorem]{Proposition}
\newtheorem{corollary}[theorem]{Corollary}

\theoremstyle{definition}
\newtheorem{definition}[theorem]{Definition}
\newtheorem{example}[theorem]{Example}
\theoremstyle{remark}
\newtheorem*{remark}{Remark}

\newcommand{\Z}{\mathbb Z}    
\newcommand{\R}{\mathbb R}    
\newcommand{\C}{\mathbb C}    
\newcommand{\PP}{\mathbb P}   
\newcommand{\T}{\mathbb T}    
\newcommand{\N}{\mathbb N}

\newcommand{\Aff}{\operatorname{Af{}f}}

\renewcommand{\O}{\mathcal O}   

\newcommand{\suchthat}{\ : \ }
\newcommand{\<}{\langle}   
\renewcommand{\>}{\rangle} 
\newcommand{\const}{\operatorname{const}}

\newcommand{\Rg}{\mathbb{R}^g}

\newcommand{\Pic}{\operatorname{Pic}}
\newcommand{\Picd}{\operatorname{Pic}^d}
\newcommand{\Div}{\operatorname{Div}}
\newcommand{\Divd}{\operatorname{Div}^d}
\newcommand{\Rtr}{\mathbb T}
\newcommand{\length}{\operatorname{length}}

\newcommand{\Sym}{\operatorname{Sym}}
\newcommand{\Val}{\operatorname{Val}}
\newcommand{\val}{\operatorname{val}}

\newcommand{\dist}{\operatorname{dist}}
\newcommand{\Log}{\operatorname{Log}}
\newcommand{\am}{{\mathcal{A}}}
\newcommand{\ignore}[1]{\relax}
\newcommand{\dd}{\partial}
\newcommand{\supp}{\operatorname{supp}}

\begin{document}

\title[Jacobians of tropical curves]
{Tropical curves, their jacobians
and Theta functions}

\author{Grigory Mikhalkin}
\address{Mathematics Department \\ University of Toronto\\ Toronto, ON M5S 2E4\\ Canada}
\email{mikha@math.toronto.edu}
\author{Ilia Zharkov}
\address{Max Planck Institut f\"ur Mathematik\\ Vivatsgasse, 7\\ D-53111 Bonn \\ Germany}
\email{zharkov@mpim-bonn.mpg.de}

\thanks{The first author is supported in part by NSERC and the second author was supported in part by NSF grant DMS-0405939}

\begin{abstract}
We study Jacobian varieties for tropical curves. These are real tori equipped with integral affine structure and symmetric bilinear form. We define tropical counterpart of the theta function and establish tropical versions of the Abel-Jacobi, Riemann-Roch and Riemann theta divisor theorems.
\end{abstract}
\maketitle

\section{Introduction}
In this paper we study algebraic curves defined
over the so-called tropical semifield $$\T=\R\cup\{-\infty\}$$
equipped with the following {\em tropical arithmetic operations} for $x,y\in\T$.
$$``x+y"=\max\{x,y\},\ ``xy"=x+y.$$
Here we use the quotation marks to distinguish between
the classical and tropical operations.
It is easy to check that each of these operations
is commutative and associative and that together
they satisfy to the distribution law.

These arithmetic operations do not turn $\T$ into a field. The
idempotency of addition, $``x+x"=x$ makes subtraction impossible.
However it does admit a division: $``\frac{x}{y}"=x-y$ for
$y\neq-\infty$. Furthermore, $\T$ has the additive zero
$``0_\T"=-\infty$\ and the multiplicative unit $``1_\T"=0$.
We say that $\T$ is a {\em semifield} as we have all the operations
except for subtraction.

It is possible to define {\em algebraic varieties over $\T$}
in a similar way as one defines algebraic varieties over
fields, e.g. $\C$ or $\R$, see \cite{Mibook}.
In this paper, we are concerned only with tropical curves
which we define in the next section.

Tropical arithmetic operations are related to classical ones
via the so-called {\em dequantization}.
Namely let us define a family of addition operations
on $\R\cup\{-\infty\}$ by
$$x\oplus_t y=\log_t(t^x+t^y),$$
where $t$ is the parameter ranging from $e$ to $+\infty$.
For any finite value of $t$ this operation is (by its definition)
induced from $\R_{\ge 0}=\{z\in\R\ |\ z>0\}$.
Clearly, $\log_t(t^xt^y)=x+y$ and thus $\R\cup\{-\infty\}$
equipped with $\oplus_t$ for addition and $+$ for multiplication
is isomorphic (as a semifield) to $\R_{\ge 0}$ with the usual
arithmetical operations. However, we have
$$\lim\limits_{t\to +\infty}\log_t(t^x+t^y)=\max\{x,y\},$$
which exhibits $\T$ as the limit semifield.

This degeneration of addition in $\R_{\ge 0}$ underlies
the following collapse of complex algebraic varieties
(see \cite{KS} and \cite{GW} for the case of Calabi-Yau varieties)
in $(\C^*)^n$).
Let
$$K=\{a(t)=\sum_{j\in I} \alpha_j t^j\ |\ \alpha_j\in\C\},$$
here $I\subset\R$ is countable and well-ordered and $a(t)$ converges
for $t\in [0,t_0]$, $t_0>0$.
Such $K$ is an example of the so-called non-Archimedean
field: there exist a valuation, i.e. a function
$$\val:K=K\to\T$$
such that $\val^{-1}(-\infty)=0_K$ and for any
$z,w\in K$ we have $\val(z+w)\le ``\val(z)+\val(w)"=\max\{\val(z),\val(w)\}$
and $\val(zw)=``\val(z)\val(w)"=\val(z)+\val(w)$.
We have the valuation map
$$\Val:(K^*)^n\to\R^n,$$
where $K^*=K\setminus\{0\}$, defined by
$\Val(z_1,\dots,z_n)=(\val(z_1),\dots,\val(z_n))$.

Let $V\subset(K^*)^n$ be an algebraic variety over $K$.
Following Kapranov \cite{Ka}
one may associate to $V$ its {\em non-Archimedean amoeba},
i.e. the image $\am=\Val(V)\subset\R^n$.
For small values of $t>0$ the variety $V$ determines
a complex algebraic variety $V_t\subset (\C^*)^n$ by
plugging $t$ to all converging series in $t$ that form
the coordinates of $V\subset(K^*)^n$.
According to Gelfand, Kapranov and Zelevinski \cite{GKZ}
one may also associate to $V_t$ its amoeba by
taking the image $\am_t=\Log_t(V_t)\subset\R^n$, where
$$\Log_t:(\C^*)^n\to\R^n$$
is defined by
$\Log_t(z_1,\dots,z_n)=(\log_t(z_1),\dots,\log_t(z_n))$.
Note that the base $t$ of the logarithm coincides
with the parameter of deformation of $V_t$.
We have the following diagram
\begin{center}
\mbox{}\xymatrix{V_t \ar[r]^{\subset} \ar[d]^{\Log_t} &
(\C^*)^n \ar[d]^{\Log_t}\\
 \am_t \ar[r]^{\subset}  & \R^n}
\end{center}
The limit of $\Log_t(V_t)\subset\R^n$ coincides with the
non-Archimedean amoeba $\am\subset\R^n$ and is an example
of tropical variety by degeneration of complex varieties.

It was shown by Kapranov \cite{Ka} that if $V\subset(K^*)^n$
is a hypersurface then $\am=\Val(V)$ depends only on the valuation
of coefficients of the polynomial defining $V$. These valuations
are elements of $\T$ and can be considered as coefficients of
the {\em tropical polynomial} defining $\am$.

The construction above gives a restricted class of tropical
varieties. In particular, they are embedded to $\R^n=(\T^*)^n$. As
the subject of this paper is tropical curves we won't give
a general definition of higher-dimensional tropical varieties
here (see \cite{Mi-ICM}). On the other hand, in this paper we define tropical
curves intrinsically and study their inner geometry that does
not depend on a particular embedding to an ambient space.

To conclude the introduction we give the definition of a particularly
useful higher-dimensional tropical variety.
\begin{definition}
The {\em tropical projective space} $\T\PP^n$ consists of the classes
of $(n+1)$-tuples of tropical numbers such
that not all of them are equal to $-\infty$
with respect to the following equivalence relation.
We say that $$(x_0:\dots:x_n)\sim (y_0:\dots:y_n),$$
if there exists $\lambda\in\T^*$ such that $x_j=y_j+\lambda$, $j=0,\dots,n$.
\end{definition}

Note that $\T\PP^n$ contains $\R^n=(\T^*)^n$:
a point $(x_1,\dots,x_n)$ corresponds to $(0:x_1:\dots:x_n)$.
Thus we have an embedding
\begin{equation}\label{Rn}
\iota_n:\R^n\subset\T\PP^n.
\end{equation}
Also we have $n+1$ affine charts
$\T^n\to\T\PP^n,$ given by the (tropical) ratio with the $j$th coordinate.

\section{Some tropical algebra}

\subsection{Tropical modules}
Recall that a commutative semigroup with zero is a set $V$ equipped with
an arithmetic operation (called addition and denoted with ``+") that is commutative,
associative and such that there exists a neutral element (called zero in $V$
and denoted with $-\infty$).
\begin{definition}
A tropical module or a {\em $\T$-module} $V$ is a commutative
semigroup with zero equipped with a map $\T\times V\to V$ (called
multiplication by a scalar) such that
\begin{itemize}
\item $``c(v+w)"=``cv+cw"$ for any $c\in\T$;
\item $``c(d(v))"=``(cd)v"$ for any $c,d\in\T$, $v\in V$;
\item $``0v"=``v"$ for the multiplicative unit $0\in\T$ and any $v\in V$;
\item if $``cv"=``dv"$ for some $c,d\in\T$, $v\in V$ then either $c=d$ or $v=-\infty$.
\end{itemize}
\end{definition}

A basic example of a tropical module is the free $n$-dimensional
tropical module $\T^n$. The addition and multiplication by a scalar
are the corresponding coordinate-wise tropical operations.
However, unlike the case with honest vector spaces over fields,
there do exist other tropical modules.

\begin{example}\label{Vn}
Let $V_n\subset\T^n$ be the submodule generated by
elements
$$e_1=(-\infty,0,\dots,0),\  e_2=(0,-\infty,0,\dots,0),\ \dots,\ e_n=(0,\dots,0,-\infty).$$
Elements of $V_n$ are linear combinations
$``\sum_{j=1}^nc_je_j"$ for some $c_j\in\T$.
\end{example}

It is easy to see that (at least as a topological space with
topology induced from $\T^n=[-\infty,+\infty)^n$), the module $V_n$
is 2-dimensional yet it is distinct from $\T^2$ (or $\T^k$ with any
other $k$). Indeed, suppose that $a_j\in\T$ and $a_1\ge a_2\ge...\ge a_n$.
Then
$$``\sum a_je_j"=(a_2,a_1,\dots,a_1).$$

For algebraic considerations we need an intrinsic definition
of dimension.
\begin{definition}\label{algdim}
Let $V$ be a tropical module. We say that its dimension is
{\em smaller than $k$} if for any elements $v_1,\dots,v_k\in V$
and any linear combination
$$v=\sum_{j=1}^k c_jv_j,$$
$c_j\in\T$, the element $v$ can be presented
as a tropical linear combination of a proper subset of
$\{v_1,\dots,v_k\}$.

We set {\em the dimension} to be equal to the maximal value of $k$
such that the dimension is not smaller than $k$.
\end{definition}

\begin{proposition}
The dimension of $V_n$ (see Example \ref{Vn}) equals 2 regardless of $n$.
\end{proposition}
\begin{proof}
Suppose that $v=\sum_{j=1}^k c_jv_j.$
We saw that each $c_jv_j$ is parameterized by two numbers $a^{(j)}_1\ge a^{(j)}_2\in\T$.
In the collection $\{c_jv_j\}$ we keep the vector $v_j$ with the maximal $a^{(j)}_1$
and another vector $v_{j'}$ with $c_{j'} v_{j'}$ having the maximal coordinate at the position of $a^{(j)}_2$.
\ignore{
Elements of $V_n$ are such $n$-tuples $v=(c_1,\dots,c_n)$ that
for some $j=1,\dots,n$
we have $c_k\ge c_j$ and $c_k=c_l$ as long as $k\neq j\neq l$.
Thus, once such $j$ is fixed we may parameterize such an element
with two parameters $t=\min\limits_{l=1,\dots,n}\{c_l\}$
and $s=\max\limits_{l=1,\dots,n}\{c_l\}$. We say that $v$ sits
on the $j$th ray at the distance $(s-t)$ from the center.

Note that if $v\in V_n$ belongs to
the $j$th ray and $v'\in V_n$ belongs to the $j'$th ray with $j\neq j'$
then $``cv+c'v'"$ for different $c,c'\in\T$ runs over
all elements sitting on the $j$th ray and distanced from the
center not greater than $v$ and all elements sitting on
the $j'$th ray and distanced from the
center not greater than $v'$.
Similarly, if $v$ and $v'$ belong to the same $j$th ray then
$``cv+c'v'"$ runs over all elements from the same ray whose
distance from the center is between that of $v$ and $v'$.
Thus if a linear combination $v=``\sum_{j=1}^k c_jv_j"$
belongs to the $j$th ray then we may re-represent it as a
tropical linear combination of only two elements from $v_1,\dots,v_k$.
}
\end{proof}
Let us introduce another useful characteristic of a $\T$-module $V$.
The {\em rank} of $V$ is the minimal number of generators of $V$.
Note that even though the modules $V_n$ have the same dimension,
they have different ranks.

\begin{proposition}
The dimension from Definition \ref{algdim} coincides with the
topological dimension of $V$.
\end{proposition}
\begin{proof}
Fixing $v_1,\dots,v_k$ gives us a piecewise-linear continuous map
$\T^k\to V$:
$$(c_1,\dots,c_n)\mapsto``\sum_{j=1}^k c_j v_j".$$
Suppose it is a surjection to $U\subset V$.
If the topological dimension of its image is $d$ then $k\ge d$.
On the other hand by the convexity argument $U$ has to be contained
in the image of the union of the $d$-dimensional coordinate subspaces
of $\T^k$.
\end{proof}

\subsection{Projectivization of a tropical module}

We may generalize the construction of $\T\PP^n$
by projectivizing {\em any} tropical module $V$, not
necessarily a free module $\T^{n+1}$.
\begin{definition}
Let $V$ be a tropical module.
Its {\em projectivization} $\PP(V)$ consists of the classes
of the non-zero elements of $v$
with respect to the following equivalence relation.
We say that $v\sim v',$
if there exists $\lambda\in\T^*$ such that $v=``\lambda v''$.
\end{definition}
Clearly we have $\PP(\T^{n+1})=\T\PP^n$.

\begin{example}
\label{Gamma-n}
The projectivization $\PP(V_n)$ of the module from Example \ref{Vn}
is called $\Gamma_n$. It will serve as our local model
for a tropical curve.
Note that by definition we have
$\Gamma_n\subset\T\PP^{n-1}$.
\end{example}

\section{Tropical curves}
\subsection{Definitions}
Any affine-linear map $A:\R^k\to\R^{k'}$ is a composition
of a linear map $L_A:\R^k\to\R^{k'}$ and a translation.
We call $A$ a {\em $\Z$-affine map} if its linear part $L_A$
is defined over $\Z$ (i.e. given by a matrix with integer entries)
no matter what is its translational part.

Let $C$ be a connected topological space
homeomorphic to a locally finite 1-dimensional simplicial complex.

\begin{definition}\label{deftropcurv}
A complete tropical structure on $C$ is
an open covering $\{U_\alpha\}$ of $C$
together with embeddings $$\phi_\alpha:U_\alpha\to\T\PP^{k_\alpha-1}$$
called {\em the charts}, $k_\alpha\in\N$,
subject to the following conditions that must hold for any $\alpha$ and $\beta$.
\begin{itemize}
\item We have $\phi_\alpha(U_\alpha)\subset\Gamma_{k_\alpha}$.
\item If $U'\subset U_\alpha$ is an open subset then
$\phi_\alpha(U')$ is open in $\Gamma_{k_\alpha}\subset\T\PP^{k_\alpha-1}$.
\item The ``finite part"
$\iota^{-1}_{k_\alpha-1}\circ\phi_\alpha\circ\phi_\beta^{-1}\circ\iota_{k_\beta-1}$
of the corresponding overlap maps
is a restriction of a $\Z$-affine linear map $\R^{k_\beta-1}\to\R^{k_\alpha-1}$.
Here we consider the map
$\iota^{-1}_{k_\alpha-1}\circ\phi_\alpha\circ\phi_\beta^{-1}\circ\iota_{k_\beta-1}$
only where it is defined. The embedding $\iota_n:\R^n\to\T\PP^n$
is taken from \eqref{Rn}.
\item If $S\subset U_{\alpha}$ is a closed set in $C$ then
$f(S)\cap \iota_{k_\alpha-1}(\R^{k_\alpha-1})$ is a closed
set in $\iota_{k_\alpha-1}(\R^{k_\alpha-1})$.
\end{itemize}
\end{definition}

\begin{figure}[htbp]
   \includegraphics{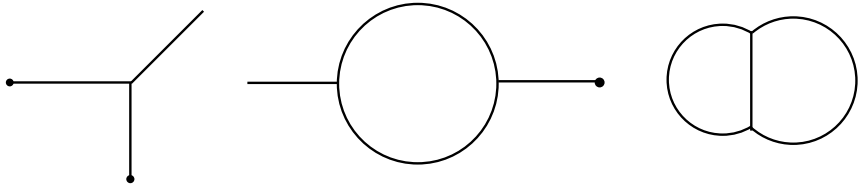}
   \caption{Examples of tropical curves of genus 0, 1 and 2.}
   \label{fig:curves}
\end{figure}

The space $C$ equipped with a complete tropical
structure is called a {\em tropical curve}.
We are especially interested in compact curves.
E.g. the first two curves on Figure \ref{fig:curves}
are not compact as each of them has an open end.
Meanwhile, it is easy to compactify these curves by
adding a point at infinity for each open end.

Clearly, the valence of a vertex of a 1-dimensional simplicial
complex is a topological invariant. Furthermore, a point inside
of an edge can be prescribed to be of valence 2 as it is possible
to introduce a new 2-valent vertex at such point.

\begin{definition}
The {\em finite part} $C^\circ$ of $C$ is the complement of
all 1-valent vertices in $C$. The 1-valent vertices are called
the {\em endpoints} of $C$.
\end{definition}

\begin{proposition}
For any chart $\phi_\alpha:U_\alpha\to\T\PP^{k_\alpha-1}$
we have
$$\phi_\alpha^{-1}(\iota_{k_\alpha-1}(\R^n))=C^\circ\cap U_\alpha.$$
\end{proposition}

This proposition follows from the last condition
in Definition \ref{deftropcurv}, the one that is
responsible for the completeness of the tropical structure.
In other words, the 1-valent vertices must sit on the boundary part
of $\T^N$.

\begin{definition}\label{tropmap}
Suppose that $C_1$ and $C_2$ are tropical curves.
A continuous map $$f:C_1\to C_2$$ is called tropical
if for every $x\in C_1$ there exists
an open neighborhood $U\ni x$ and charts
$U_\alpha\supset U$ and $U_\beta\ni f(x)$,
$U_\alpha\subset C_1$, $U_\beta\subset C_2$,
such that
$\iota_{k_\beta}\circ\phi_{\beta}\circ f\circ \phi^{-1}_\alpha\circ\iota_{k_\alpha}$
is a restriction
of a $\Z$-affine map $\R^{k_\alpha}\to\R^{k_\beta}$.

Tropical curves are called isomorphic if there exist
mutually inverse tropical maps between them.
\end{definition}

\subsection{Equivalence of tropical curves}

It turns out that some non-isomorphic tropical curves
serve as two distinct models for the same variety.
It is convenient to introduce an equivalence relation
identifying such models.

Let $C$ be a tropical curve, $x\in C^\circ$ be a point
in the finite part of $C$.
Restricting to another chart if needed
we can find a chart $\phi_\alpha:U_\alpha\to\T\PP^{k_\alpha-1}$
such that $\phi_\alpha(x)=(1_{\T}:\dots:1_{\T})=(0:\dots:0)$
is the ``center" of the curve $\Gamma_{k_\alpha}\subset\T\PP^{k_\alpha-1}$.

The projection along the last coordinate defines a continuous map
$\lambda_{k_\alpha}:\T\PP^{k_\alpha}\setminus\{(-\infty:\dots:-\infty:0)\}
\to \T\PP^{k_\alpha-1}.$
The map
\begin{equation}\label{lambdagamma}
\lambda_{\Gamma_{k_\alpha}}:
\Gamma_{k_\alpha+1}\to\Gamma_{k_\alpha}
\end{equation}
obtained by restricting $\lambda_{k_\alpha}$ to $\Gamma_{k_\alpha+1}$
and extending it to $(-\infty:\dots:-\infty:0)$ is continuous.

The map $\lambda_{\Gamma_{k_\alpha}}$
contracts the interval $[(0:\dots:0),(-\infty:\dots:-\infty:0)]$
to $(0:\dots:0)$ and is a homeomorphism to its image when restricted
to $\Gamma_{k_\alpha+1}\setminus ((0:\dots:0),(-\infty:\dots:-\infty:0)]$.
Let $$\tilde{U_\alpha}=\lambda_{\Gamma_{k_\alpha}}^{-1}(\phi_\alpha(U_\alpha)).$$

We form $\tilde{C}=C_x$ by replacing the neighborhood $\tilde{U_\alpha}$ with
$U_\alpha$. Clearly, the tropical structure naturally extends
to the whole curve $\tilde{C}$. Topologically, $\tilde{C}$ is
the result of gluing of a closed interval $[0,1]$
to $C$ by identifying $0\in [0,1]$ with $x\in C$. Denote with
$\tilde{x}\in\tilde{C}$ the point corresponding to $1\in [0,1]$.
(See \cite{Mi-ICM} for a more general procedure that can be used
to define higher-dimensional tropical varieties.)

\begin{definition}
The map
\begin{equation}\label{contract}
\tilde{C}\to C
\end{equation}
induced by \eqref{lambdagamma}
is called the {\em elementary
equivalence} of tropical curves. Two curves $C_1$ and $C_2$
are called {\em tropically equivalent} if they can be connected
with a sequence of elementary equivalences (in any direction).
\end{definition}

This equivalence relation provides a convenient tool
for working with curves with marked points.
If $x\in C^\circ$ is in the affine part then
we may replace $C$ with $\tilde{C}=C_x$ and $x$ with
$\tilde{x}$. In this way we may replace a curve $C$ with
a collection of marked points on it with an equivalent curve $\bar{C}$,
so that all the marked points will be the endpoints of $\bar{C}$.
This trick turns out to be very useful for the setup of the Gromov-Witten
theory, cf. \cite{Mibook}.

Here we would just like to note that putting marked points
to the end of the equivalent curve can be used to puncture
a tropical curve. If we need to remove a point $x\in C^\circ$
then first we replace $(C,x)$ with an equivalent pair $(\tilde{C},\tilde{x})$
and then treat $\tilde{C}\setminus\tilde{x}$ as the result of puncturing
$C$ at $x$. Note that $\tilde{C}\setminus\tilde{x}$ is a tropical curve
while $C\setminus x$ itself is not.
Furthermore, the contraction map \eqref{contract} can be considered
as a morphism analogue of a Zariski-open set (in the Grothendieck topology style), see
\cite{Mibook} for details.

\subsection{Tropical curves as metric graphs}

A tropical curve turns out to carry the same information as
the so-called {\em metric graph}. Recall that the edges adjacent to
one-valent vertices of a graph are called {\em leaves}.
Other (non-leaf) edges are called {\em inner}.
A metric graph (cf. e.g. \cite{Biomath}) is a graph
equipped with finite positive length for all its inner edges.
The leaves of a metric graph are prescribed the (positive) infinite length.

Let $\Gamma$ be a connected finite graph and ${\mathcal V}_1(\Gamma)$
be its 1-valent vertices (or {\em endpoints}).
We may consider
$$\Gamma^\circ=\Gamma\setminus{\mathcal V_1}(\Gamma).$$
If $\Gamma$ is a metric graph then $\Gamma^\circ$ is a complete
metric space (with an inner metric). Vice versa, a complete
metric space homeomorphic to $\Gamma^\circ$ for some finite graph $\Gamma$
defines a metric graph.

Recall that if $X$ is any topological space one can consider
its maximal compactification $\bar{X}$ by adding a point ``at infinity"
for each end of $X$. Clearly, $\Gamma$ and
$\bar{\Gamma^\circ}$ are homeomorphic if $\Gamma$ is a finite graph.

\begin{proposition}
There is a natural 1-1 correspondence between compact tropical curves
and metric graphs.
\end{proposition}
\begin{proof}
Let $L\subset\R^k$ be a line of a rational slope
and let $\xi\in\Z^k$ be a primitive (i.e. not divisible by
an integer)
vector parallel to $L$. Define the distance between
points $x,y\in L$ to be $\frac{||x-y||}{||\xi||}$,
where $||\ .\ ||$ stands for any norm in the vector space $\R^k$.
This metric is preserved by the overlap maps
$\phi_\alpha\circ\phi_\beta^{-1}
|_{\phi_\beta(U_\alpha\cap U_\beta)\cap\R^{k_\beta}}$.
This gives us a natural inner metric on a tropical curve.

Vice versa, if $C$ is a metric graph then we set
the metric-preserving charts to $\R$ for the interior of the edges;
charts to $\Rtr$ for small neighborhoods of 1-valent vertices and
isometric embedding charts to $\Gamma_k$ for small
neighborhoods of multivalent vertices.
\end{proof}

\subsection{Regular and rational functions on curves}

By the very definition
a tropical curve $C$ is a topological space homeomorphic to a graph
equipped with an additional geometric structure which we call
{\em the tropical structure}.
It allows one to define a sheaf of regular function on $C$.
We say that a function $\R^n\to\R\subset\T$ is {\em $\Z$ affine-linear}
if it is obtained from a function $\R^n\to\R$ that is linear over $\Z$
by adding a real constant.

Let $U\subset C$ be an open set.
\begin{definition} A function $f:U\to\T$ is called {\em regular}
if for any $x\in U$ there exists a chart
$$\phi_\alpha:U_\alpha\to\T^n\subset\T\PP^n,$$
$n=k_\alpha-1$, $U_\alpha\ni x$, a tropical polynomial
$g:\T^n\to\T$,
$$g(x_1,\dots,x_n)=
``\sum_{(j_1,\dots,j_n)\in \mathcal A} a_{j_1\dots j_n}x_1^{j_1}\dots x_n^{j_n}"=
\max\limits_{(j_1,\dots,j_n)\in \mathcal A} \{a_{j_1\dots j_n}+j_1 x_1+\dots+ j_n x^n\},$$
where $\mathcal A$ is a finite subset of $\Z^n_{\geq 0}$,
such that the function
$f\circ\phi^{-1}_\alpha - g$ (which is defined on $\phi_\alpha(U)$)
restricted to $\R^n\subset\T^n$ is $\Z$ affine-linear.
Here $\T^n\subset\T\PP^n$ is embedded given by
$(x_1,\dots,x_n)\mapsto (0:x_1:\dots:x_n)$.
If the function $g$ may be chosen to be a constant then
$f$ is called {\em an affine-linear function on $U$}.

A function $h:U\to \T$ is called {\em rational} if there exist two
regular functions $f_1,f_2:U\to \T$ such that
$$h=``\frac{f_1(x)}{f_2(x)}"=f_1(x)-f_2(x)$$
for any $x\in U\cap C^\circ$.
\end{definition}

Clearly, regular functions form a sheaf of tropical algebras on $C$.
The same holds for rational functions.
The sheaf of regular functions is called the {\em structure sheaf}
and is denoted with $\O^{\T}_C$.

\subsection{Projective tropical curves, tropical curves in $\R^n$
and in other higher-dimensional tropical varieties}
Linear morphisms maybe used to embed or immerse
tropical curves into {\em tropical toric varieties}
that arise as compactifications of $\R^n=(\T^{\times})^n$.
Our main example of such compactification is $\T\PP^n$.

Let $f:C\to\T$ be an affine-linear function.
For every edge $E\subset C$ we may define {\em the slope}
of $f$ once we fix an orientation of $E$. Indeed if $\xi$
is a primitive tangent vector consistent with the choice
of the orientation then the slope is just the partial
derivative $\frac{\dd f}{\dd\xi}$. If we do not fix the
orientation of $E$ then the slope is only defined up to
the sign.

\begin{definition}
A map $C\to \T^n$ is called a {\em linear morphism}
if it is given in coordinates by affine-linear functions.
A map $C\to \T\PP^n$ is called a {linear morphism}
if locally (with respect to the charts $\T^n\subset\T\PP^n$)
it is given by linear morphisms.
\end{definition}
\begin{remark}
This definition agrees with a more general definition given
in \cite{Mi-ICM}. It is easy to show that every map
$C\to\T\PP^n$ locally given in coordinates by regular functions
admits a resolution by a linear morphism $\tilde{C}\to\T^n$
for an equivalent tropical curve $\tilde{C}$ with a contraction
$\tilde{C}\to C$.
\end{remark}

Let $C$ be a compact tropical curve and
$h:C\to\T\PP^n$ be a linear morphism. The image $h(E)$ of an edge
$E\subset C$ is either a point (if the slope of all coordinates is
zero) or a straight interval with a rational slope (as the slopes
of all coordinates are integers). We define {\em the weight} $w(E)$
as the GCD of the coordinate slopes of $E$.
The collection of the slopes of all coordinates gives
the {\em weight vector} $\xi_E\in\Z^n$ (defined up to
a sign unless we specify the orientation of $E$).
The ratio $\frac{\xi_E}{w(E)}$ is a
primitive integer vector parallel to $E$.

The image $h(E)$, if it not a single point, is contained in the finite part $\R^n\subset\T\PP^n$
if and only if both vertices adjacent to $E$ have valence greater than 1.
All 1-valent vertices of $C$ are mapped to the boundary part $\T\PP^n\setminus\R^n$.
Let $P\subset C$ be a vertex and $E_1,\dots,E_k$ are the edges adjacent
to $P$. If $k>1$ then $h(P)\in \R^n$ and $\sum_{j=1}^k\xi_k=0$.

\begin{figure}[htbp]
   \includegraphics{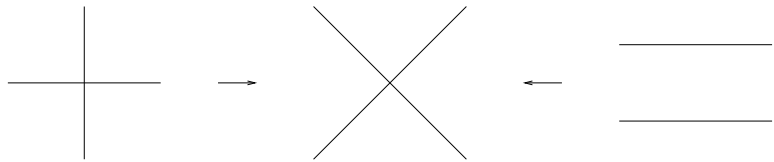}
   \caption{\label{vr2} A cross figure (in the center) obtained as the image of
   two different linear morphisms.}
\end{figure}

Recall that a {\em leaf} of a graph is an edge adjacent to a 1-valent vertex.
The {\em degree} of a projective curve is the number of leaves
(counted with some weight, cf. \cite{Mi07}) adjacent to any of the $(n+1)$
components of the boundary divisor $\T\PP^n\setminus\R^n$.
Note that linear morphisms of different curves may have the
same image, see Figure \ref{vr2} for an example.

If $X\subset\T\PP^n$ is a projective tropical variety (see \cite{Mi-ICM}) then
we define a (parameterized) tropical curve in $X$ as a linear morphism
$C\to\T\PP^n$ such that its image is contained in $X$. Note that to embed
$X$ to a projective space we often need to pass to an equivalent model of $X$.

\section{Divisors and line bundles}

\subsection{Tangent vectors and 1-forms}

Let $\xi$ be a primitive integer vector tangent
to an interval in $\phi_\beta(U_\alpha\cap U_\beta)\subset\R^{k_\beta}$.
Note that the overlap maps $\phi_\alpha\circ\phi_\beta^{-1}
|_{\phi_\beta(U_\alpha\cap U_\beta)\cap\R^{k_\beta}}$
takes $\xi$ to a primitive integer tangent vector
to an interval in $\phi_\alpha(U_\alpha\cap U_\beta)\subset\R^{k_\alpha}$
(since the overlap map admits a $\Z$-affine inverse).
Thus we have a well-defined notion of a {\em primitive tangent
vector} to a point $p\in C$.

If $p$ is a point in the interior of an edge then it has two
primitive tangent vectors $\xi_1$ and $\xi_2$. We have $\xi_2=-\xi_1$.
An {\em integer vector tangent to $C$ at $p$} is
an integer multiple of $\xi_1$ (or $\xi_2$).

If $p$ is a vertex of $C$ then we can distinguish between outward
and inward primitive tangent vectors at $p$. If $p$ is a $k$-valent
vertex then it has $k$ outward tangent vectors $\xi_1,\dots,\xi_k$
with $$\sum_{j=1}^k\xi_j=0$$ (the curve $C$ is an abstract
curve now, but this equality makes sense in any affine chart of $C$). An integer
vector tangent to $C$ at $p$ is an integer linear combination of
$\xi_1,\dots,\xi_{k-1}$. It defines an integer vector in
$\R^{k_\alpha}$ for any chart $U_\alpha\ni p$.

Note that there is no difference between a 2-valent vertex
and an interior point of an edge. Both of them have affine
charts to $\R$. At our convenience we may introduce extra
2-valent vertices by subdividing an edge or do the opposite
operation.

\begin{definition}
Let $\Aff_\R$ be the sheaf of affine functions with real (i.e., not
necessarily integral) slope. Define the {\em real} cotangent local system
${\mathcal T}^*$ on $C$ by the following exact sequence of sheaves:
\begin{equation*}
0 \longrightarrow \R  \longrightarrow \Aff_\R \longrightarrow {\mathcal T}^*\longrightarrow
0,
\end{equation*}
A {\em 1-form} on $C$ is a global section of ${\mathcal T}^*$.
\end{definition}

In particular, since a 1-form $\omega$ is locally constant it must be balanced at any $p\in C$:
\begin{equation}\label{eq:balance}
\sum_{i=1}^{\val(p)} \omega(\xi_i) = 0,
\end{equation}
where $\val(p)$ is the valence of $p$ and the $\xi_i$'s are the outward primitive integral
tangent vectors at $p$.

\begin{remark} The fact that the forms have to be constant on edges can be interpreted in the degeneration picture as follows. The
regular forms on complex curves that survive in the limit are
of the form $\frac\alpha{2\pi i} \frac{dz}{z}$ on long cylinders
(edges to be). The coefficient $\alpha$ would be the value of
the limiting form $\omega$ on the primitive vector tangent to the
edge. The balancing condition (\ref{eq:balance}) reflects the
fact that the sum of residues of a rational form on $\PP^1$ (a
vertex to be) is zero.
\end{remark}

\subsection{Divisors}
A {\em divisor} on $C$ is a formal linear combination of points in $C$ with integer coefficients: $D=\sum a_i p_i$, where $a_i\in \Z$. We define its {\em degree} in the usual way:
$\deg(D) := \sum a_i$.
The divisor $D$ is called {\em effective} if all the coefficients are non-negative.
Clearly, the degree of an effective divisor is non-negative.
The set of all divisors form an
abelian group which we will denote by $\Div(C)$.

Given a rational function
$f$ on an open subset $U\subset C$ one can consider the divisor
\begin{equation}\label{eq:pr_divisor}
(f):=\sum_{p\in U} \left(\sum_{i=1}^{k(p)}\frac{\partial f}{\partial \xi_i}(p)
\right) p.
\end{equation}
If $U=C$, i.e. $f$ is a global function, then $(f)$ called {\em
principal}. We say that two divisors are {\em linearly equivalent}: $D_1\sim
D_2$, if $D_1-D_2$ is principal.

\begin{figure}[htbp]
   \includegraphics{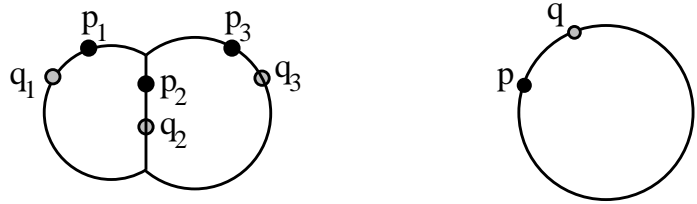}
   \caption{Linear equivalence: $p_1+p_2+p_3\sim q_1+q_2+q_3$, but $p\nsim q$.}
   \label{fig:lin_equivalence}
\end{figure}

\begin{proposition}\label{lemma:residue}
Let $U$ be an open set in a tropical curve $C$ whose boundary consists of the points $z_1,\dots,z_k$. Let $\nu_i$ be the outward (that is pointing away from $U$)
primitive tangent vectors at the $z_i$. Let $f$ be a rational function on $U$ with
divisor $(f)=\sum a_j p_j$. Then
\begin{equation}\label{eq:degree}
\sum a_j =\sum_{i=1}^k \frac{\partial f}{\partial \nu_i}(z_i).
\end{equation}
\end{proposition}
\begin{proof}
We observe that the formula is additive when connecting pieces of the curve at
points other than the $p_i$. Note also, that the formula trivially holds for
affine linear functions. On the other hand at the points $p_i$ it is just the
definition of the order $a_i=\sum_{j=1}^{\val(p_i)}\frac{\partial f}{\partial
\xi_j}(p_i)$.
\end{proof}

\begin{corollary}
A principal divisor on a compact curve has degree zero.
\end{corollary}

\subsection{Line bundles}
\begin{definition}
A {\em line bundle} on $C$ is an ${\mathcal O}^*$-torsor, where the sheaf of
affine functions ${\mathcal O}^*$ acts on the structure sheaf by tropical
multiplication. Equivalently a line bundle $L$ is a topological space together
with the following data:
\begin{itemize}
\item A continuous projection $\pi: L\to C$ with fibers $\Rtr$.
\item Every point $p\in C$ has an open neighborhood $U_p\ni p$ and families of
{\em trivialization} maps $\{\phi\}_V: \pi^{-1} (V)\cong V\times \Rtr$ for any
open $V\subset U_p$, which are the restrictions from  $\{\phi\}_{U_p}: \pi^{-1}
(U_p)\cong U_p\times \Rtr$, and such that any two trivializations
$\phi_1,\phi_2\in \{\phi\}_V$ differ by (tropical multiplication by) an affine
linear function in $V$.
\end{itemize}
\end{definition}

Given a sufficiently fine open covering $\{U_\alpha\}$ of $C$ a line bundle can
be specified by transition functions $f_{\alpha\beta}:U_\alpha\cap U_\beta \to
\R$ between local trivializations. The functions $f_{\alpha\beta}$ are affine linear
and satisfy the usual cocycle condition. Another choice of trivializations will
result in changing the cocycle by a coboundary. A continuous map $L_1\to L_2$
which respects the projection and local trivialization up to affine functions
is an isomorphism of line bundles. Thus, as in the classical geometry, the
isomorphism classes of line bundles are parameterized by the first \v{C}ech
cohomology $H^1(C,\O^*)$.  This is a group and we will refer to it as the
Picard group $\Pic(C)$.

\begin{proposition}
Equivalent tropical curves have the same Picard groups.
\end{proposition}

\subsection{Local sections, global sections and linear systems}
Given an open subset $U\subset C$ a section $s:U\to\pi^{-1} (U_p)$ is {\em
regular}, respectively {\em rational},
if for any open $V\subset U$ and
a trivialization $\phi: \pi^{-1} (V)\cong V\times \Rtr$, $s$ becomes a
regular, resp. rational, function on $V$. The notion does not depend on
trivializations. Neither does the formal sum in (\ref{eq:pr_divisor}). Thus a
global rational section of $L$ defines a divisor on $C$. There is the usual
correspondence between line bundles and divisors up to equivalence:

\begin{proposition}\label{prop:divisors-bundles}
Every divisor defines a line bundle together with its
rational section. This section is defined uniquely up
to adding a constant (i.e. up to tropical multiplication).

Conversely, every line bundle $L$ has a
rational section. The divisors of any two rational sections of $L$ are
linearly equivalent.
\end{proposition}
\begin{proof}
Let $D=\sum a_i p_i$ be a divisor. We can cover $C$ by open sets $U_i$ so that
each point $p_i$ is contained in a unique $U_i$. Then for each $U_i$ we choose
a rational function $f_i$ whose order of zero (pole) at $p_i$ is $a_i$. The
incompatibilities over the intersections define a \v{C}ech cocycle with values
in affine linear functions. Hence we get a line bundle. The rational section
is given by the collection $\{f_i\}$.

Conversely, we can choose $g$ disjoint open intervals
$N_i$ in the interiors of edges such that $C\setminus {\cup N_i}$ is a (connected) tree (cf. \S \ref{break_points}). Also choose $g$ closed intervals $N^\epsilon_i\subset N_i$.
Then a given line bundle $L$ can be
trivialized over $C\setminus{\cup N^\epsilon_i}$ and over each of the $N_i$. Choose a
rational function $f$ on $C\setminus{\cup N^\epsilon_i}$ which is affine linear on the $N_i \setminus N^\epsilon_i$. To extend $f$ to a section $s$ of $L$ over $C$ one needs integral PL functions $f_i$ on the $N_i$ equal to given affine linear functions on the pairs of ends $N_i \setminus N^\epsilon_i$. This is always possible. The divisor of $s$ is $(f)+(f_1)+\dots+(f_g)$.
\end{proof}

If a divisor $D$ has degree $d$ we will say that the corresponding line bundle $L(D)$
has degree $d$ as well. The groups of degree $d$ divisors and line bundles will
be denoted by $\Div^d(C)$ and $\Pic^d(C)$ respectively. Now we can reformulate
Proposition \ref{prop:divisors-bundles}.
\begin{corollary}
For every $d$ there is a canonical isomorphism $\Pic^d(C)\cong \Div^d(C)/\sim$.
\end{corollary}

We will now look at the space of global regular
sections of a given line bundle. Let $D_0$ be a divisor and let $s_0$ be the
corresponding rational section (unique up to an additive constant) of the
line bundle $L:=L(D_0)$. Then any other rational section $s$ of $L$ is given
by $s=s_0+f$, where $f$ is a rational function. Of course, the section $s$ is
regular if the divisor $D_0+(f)$ is effective.

The space of global sections $\Gamma(C,L)$ has the structure of a
$\Rtr$-module. Given $s_1,s_2\in \Gamma(C,L)$ one can take their tropical sum
$\max\{s_1,s_2\}$ by choosing trivializations over open sets $U\subset C$ and
considering $s_1,s_2$ as functions on $U$. The result is a section of the same
bundle and it is independent of the trivializations. Adding a constant to a
section is well defined as well.

Given a divisor $D$ we denote by $|D|$ its complete linear system, that is the collection of all {\em effective} divisors linearly equivalent to $D$. Clearly, the set $|D|$ is the tropical projectivization of the $\Rtr$-module $\Gamma(C,L)$ for the corresponding line bundle $L$. We will study linear systems in more details in the Riemann-Roch section.

\subsection{Universal covering and fundamental domains}\label{break_points}
For the proof of Jacobi Inversion and Riemann Theorem it will be convenient to work with
the universal covering $\hat C\to C$. Here is a way to fix a fundamental domain.

Let us choose $g$ points $z_1,\dots,z_g$ in $C$  equipped with outward primitive tangent vectors $\eta_1,\dots,\eta_g$.
Let $z_i^\epsilon$ be the point on $C$ which is distance $\epsilon$ from $z_i$ into the direction $\eta_i$.
For a small $\epsilon>0$ we may consider the set
$C\setminus\{z_1^\epsilon,\dots,z_g^\epsilon\}$. If it is a (connected) tree
then we refer to $\{z_i,\eta_i\}$ as a set of {\em break points} and
call $D=z_1+\dots+z_g$ the {\em break divisor}.

The whole business with $z_i^\epsilon$ is only needed to serve the case when $z_j$ are
points of valence greater than 2. In this case some of $z_j$ may coincide if the vectors
$\eta_j$ are different.

\begin{figure}[htbp]
   \includegraphics{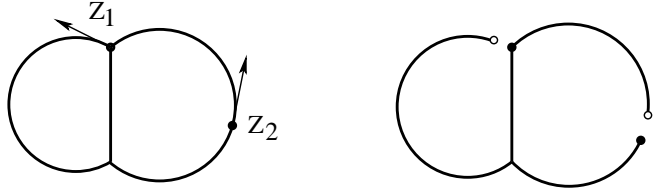}
   \caption{Break points and the fundamental domain $T$.}
   \label{fig:break}
\end{figure}

The choice of break points on $C$ is equivalent to the choice of a connected fundamental domain $T\subset \hat C$ (see Fig. \ref{fig:break}). Indeed, given $T$ the break points are the missing boundary points of $T$ with $\eta_i$ pointing inside, see Fig. \ref{fig:break}. On the other hand, given a collection of break points $z_i,\eta_i$, the fundamental domain $T$ is defined by the property that the lifting map $C\to T\subset \hat C$ is continuous except at the $z_i$ from the direction of $\eta_i$.

In general, given a (not necessarily connected) fundamental domain $T\subset \hat C$ we can define a {\em pseudo-break divisor} as follows. Using the identification of points in $T$ with those in $C$ we set
\begin{equation}\label{eq:pseudo_break}
D_T=\sum_{z\in C} (\val_C(z)-\val_T(z)) z,
\end{equation}
where $\val_C(z)$ and $\val_T(z)$ stand for the valence of $z$ in $C$ and $T$, respectively. Clearly a pseudo-break divisor $D_T$ is a break divisor if $T$ is connected (a tree). Also since $T$ contains no cycles a simple Euler characteristic count gives the following criterion.

\begin{lemma}\label{lemma:pseudo_break}
A pseudo-break divisor $D_T$ is a break divisor if it is of degree $g$.
\end{lemma}

We remark that given a break divisor $D=z_1+\dots+z_g$ there is more than one way to choose tangent vectors $\eta_i$ to make $z_i, \eta_i$ into a collection of break points. This choice will affect the choice of the fundamental domain $T$. In any case, $T$ is equal to the limit of the fundamental domains $T^\epsilon$ (induced from the break points $z_i^\epsilon$) when $\epsilon\to 0$, and the corresponding break divisors $D_{T^\epsilon}$ tend to $D_T$.

\section{Tori and polarizations}
A tropical {\em torus} $X$ is the quotient $\R^n/\Lambda$, where
$\R^n$ is considered with an integral affine structure, i.e. with
a fixed integral lattice $\Z^n\subset\R^n$, and $\Lambda\subset\R^n$
is (another) lattice. We will denote with $(\R^n)^*$ the dual space, and let $(\Z^n)^*$ and $\Lambda^*$ be the dual lattices to $\Z^n$ and $\Lambda$, respectively.

\begin{remark}
Tropical tori, especially the ones endowed with positive polarization (see below) had been considered (under different names) in the framework of degenerations of Abelian varieties (cf. \cite{Mumford_abelian_varieties}, \cite{Chai-Faltings},  \cite{Alexeev-Nakamura}, \cite{Alexeev}) long before tropical geometry came into existence.
\end{remark}

\subsection{Line bundles and polarizations} We can repeat most
of our definitions and results about line bundles and sections
for curves.  A line bundle on a torus is defined by an
element in $H^1(X,\O^*)$, where $\O^*=\Aff$ is the sheaf of affine
linear functions on $X$ with integral slope.

\begin{remark}
Note however that as in classical geometry not every line bundle on a tropical torus has a rational section.
\end{remark}

Now consider a short exact sequence of sheaves:
\begin{equation}
0 \longrightarrow \R  \longrightarrow \Aff  \longrightarrow {\mathcal T}^*_\Z \longrightarrow 0,
\end{equation}
where $\R$ is locally constant sheaf and sections of ${\mathcal T}^*_\Z \cong (\Z^n)^*$
are locally constant  integral 1-forms.
The induced map of the long exact sequence
\begin{equation}
c: H^1(X,\Aff) \longrightarrow H^1(X,{\mathcal T}^*_\Z)
\end{equation}
is called the {\em Chern class} map. A choice of a class $[c]\in H^1(X,{\mathcal T}^*_\Z)
\cong \Lambda^*\otimes (\Z^n)^*$ which is in the image of the Chern class map
is called a {\em polarization} of $X$. Using the natural isomorphism
$\Lambda^*\otimes (\Z^n)^*\cong \operatorname{Hom}(\Lambda, (\Z^n)^*)$ one can
think of $[c]$ either as a bilinear form on $\R^n$ or as a map
$[c]:\Lambda\to(\Z^n)^*$.

We claim that the set of polarization classes is formed by the maps
$\Lambda\to(\Z^n)^*$ which are {\em symmetric} as bilinear forms on
$\R^n$. To see this we analyze the coboundary map in the induced
long exact sequence:
\begin{equation}
\delta: H^1(X,{\mathcal T}^*_\Z) \longrightarrow H^2(X,\R).
\end{equation}
By identifying $H^1(X,{\mathcal T}^*_\Z) \cong \Lambda^*\otimes (\Z^n)^*$ and
$H^2(X,\R)\cong \wedge^2 (\R^n)^*$ one can deduce that $\delta$ is the
restriction of the skew-symmetrization map $(\R^n)^*\otimes (\R^n)^* \to
\wedge^2 (\R^n)^*$.

\begin{remark} For a general affine manifold the coboundary map $\delta$ is
given by a wedge product with some characteristic class $\rho$ which takes
values in the cohomology with coefficients in the tangent local system $\mathcal T$.
This  so called {\em radiance obstruction class} has the meaning of the
translational part of the monodromy representation (cf. \cite{KS04}, Section
2.2). In our case the local system $\mathcal T$ is trivial because the linear part of
the monodromy on the torus is trivial. Thus, $H^1(X,\mathcal T)$ can be canonically
identified with $(\R^n)^*\otimes \R^n$, and $\rho$ is the identity element
there.
\end{remark}

We can also describe the set of line bundles with a given Chern class. Notice
that the map $H^0(X,\R)\rightarrow H^0(X,\Aff)$ is an isomorphism. Hence
$H^0(X,{\mathcal T}^*_\Z)\rightarrow H^1(X,\R)$ is the natural inclusion
$(\Z^n)^*\hookrightarrow (\R^n)^*$. Thus the Picard group of $X$ with a fixed
Chern class can be identified with the dual torus $(\R^n)^*/(\Z^n)^*$.

If the quadratic form defined by $[c]$ is positive definite we call
$\R^n/\Lambda$ a {\em (polarized) tropical Abelian variety}
(cf. \cite{Alexeev-Nakamura}, \cite{Alexeev}). The {\em index} of
the polarization $[c]$ is the index of the image $[c](\Lambda)$ in $(\Z^n)^*$.
We call a polarization {\em principal} if it has index one, that is, if the map
$[c]:\Lambda\to(\Z^n)^*$ is an isomorphism.

\begin{theorem}\label{thm:polarization}
Let $L$ be a line bundle on $X=\R^n/\Lambda$ which defines a principal
polarization. Then the space of sections $\Gamma(X,L)$ is one dimensional.
\end{theorem}
\begin{proof}
Suppose $L$ has the Chern class $[c]:\Lambda\to(\Z^n)^*$. A global section
$\Theta$ of $L$ can be viewed as a {\em convex} PL function on $\R^n$ subject
to some quasi-periodicity condition:
$$\Theta(x+\mu)=\Theta(x) + \<[c](\mu), x\> + \beta(\mu), \quad \text{for any }
\mu\in \Lambda.$$
By applying quasi-periodicity twice we can see that $\beta$ has to be quadratic.
Namely,
$$\beta(\mu_1+\mu_2)=\beta(\mu_1)+\beta(\mu_2)+[c](\mu_1,\mu_2).$$
That is, $\beta(\lambda)=\frac12 [c](\lambda,\lambda)$ + linear part. The linear
part is responsible for shifts of the line bundle.

On the other hand, as a convex function $\Theta(x)$ is completely determined by
its Legendre transform $\hat\Theta:(\Z^n)^*\to\R$. Now we show that a value of
$\hat\Theta(0)$ determines $\hat\Theta(m)$ at any $m\in \operatorname{Im}[c]
= (\Z^n)^*$. Say $m=[c](\mu)$ for some $\mu\in\Lambda$, then
\begin{multline*}\hat\Theta(-m)=\max_{x\in\R^n}\{\<-[c](\mu), x\>-\Theta(x)\}
=\max_{x\in\R^n}\{-\Theta(x)+\beta(\mu)\}= \hat\Theta(0)+\beta(\mu).
\end{multline*}
Performing the inverse of the Legendre transform
gives
\begin{multline*}
\Theta(x):=\max_{m\in(\Z^n)^*}\{\<m,x\>-\hat\Theta(m)\}
=\max_{\lambda\in\Lambda}\{\<-[c](\lambda), x\>-\hat\Theta(0)-\beta(-\lambda)\}.
\end{multline*}
This is the unique section of $L$ determined upto a additive constant by
$\hat\Theta(0)$.
\end{proof}

\begin{remark}
For a general tropical Abelian variety, a section of a positive line bundle $L$,
thought of the quasi-periodic PL function $\Theta$ on $\R^n$, is completely
determined by specifying values of its Legendre transform $\hat\Theta$ on every
representative of the quotient $(\Z^n)^*/[c](\Lambda)$. Thus, in general, the
dimension of $H^0(X,L)$ is given by the index of the polarization.
\end{remark}

\subsection{Theta functions}
The unique (up to constant) section from Theorem \ref{thm:polarization}
will be the tropical analog of the classical Riemann's {\em theta function}.
Let $\Lambda$ be a lattice in $\R^n$ and $Q: \R^n\otimes \R^n \to \R$
be a positive definite symmetric bilinear form.
As implicitly suggested in \cite{Alexeev-Nakamura} we define
$$\Theta(x):=\max_{\lambda\in\Lambda} \{Q(\lambda,x) - \frac12
Q(\lambda,\lambda)\}, \quad x\in \R^n.$$
The maximum here always exists since $Q$ is positive definite. From the
definition we can readily see that $\Theta(x)$ is an even function:
$\Theta(-x)=\Theta(x)$. It also satisfies the following functional equation:

\begin{lemma} \label{theta_period}
$\Theta(x+\mu)=\Theta(x) + Q(\mu,x) + \frac12 Q(\mu,\mu)$, for any $\mu\in
\Lambda$.
\end{lemma}
\begin{proof}
This calculation is in a sense ``Legendre dual'' to the the proof of Theorem
\ref{thm:polarization}. By considering the effect of translation on each term
we have:
\begin{multline*}
Q(\lambda,x+\mu) -\frac12 Q(\lambda,\lambda)\\
=Q(\lambda-\mu,x) +Q(\mu,x)+Q(\lambda,\mu) -\frac12 Q(\lambda-\mu,\lambda-\mu)
-Q(\mu,\lambda) +\frac12 Q(\mu,\mu)\\
=Q(\lambda-\mu,x) - \frac12 Q(\lambda-\mu,\lambda-\mu) +Q(\mu,x)+\frac12 Q(\mu,\mu).
\end{multline*}
Relabeling the terms and combining them into the tropical sum completes the
proof.
\end{proof}

The theta function above can be defined on $\R^n$ for an arbitrary
positive symmetric bilinear form $Q$. But it is {\em regular
(holomorphic)} in tropical sense
only if the form $Q$ is integral in the sense that the image of the
induced map $\tilde Q:\Lambda\to(\R^n)^*$ ends up being inside the
lattice $(\Z^n)^*$. From now on we restrict our attention to the
integral forms $Q$ only. This is when the theta function defines a
polarization on the tropical torus $\R^n/\Lambda$. Moreover, we will
be primarily interested in principal polarizations.

\begin{figure}[htbp]
   \includegraphics{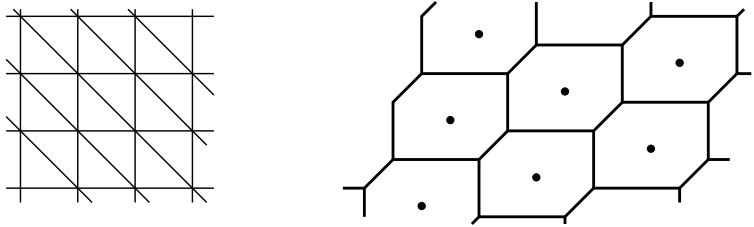}
   \caption{$\hat\Theta(m)$-induced decomposition of $(\R^n)^*$ and a
corresponding $\Lambda$-periodic theta divisor.}
   \label{fig:theta}
\end{figure}

A tropical hypersurface in $\R^n$ defined by the corner locus of the theta
function gives a honeycomb-like periodic cell decomposition of $\R^n$, the so-called Voronoi decomposition (cf. Fig. \ref{fig:theta} and \cite{Alexeev-Nakamura}). It is dual to the decomposition of $(\R^n)^*$ with vertices in $(\Z^n)^*$ induced by the
Legendre transform $\hat\Theta(m)$, the so-called Delauney decomposition.
The principal polarizations can be geometrically characterized by that
the volume of each maximal Voronoi cell is equal to the affine volume of the torus
$\R^n/\Lambda$, which, in turn, is equal to $\det \tilde Q$.

The theta function can be considered as a section of the polarization line
bundle $L$ on $X$ whose corresponding \v{C}ech 1-cocycle is defined by the
``automorphy factors'' $Q(\mu,x) + \frac12 Q(\mu,\mu)$ in the functional
equation. Locally $\Theta$ is a tropical polynomial in $\R^g$ and thus
defines a tropical hypersurface $[\tilde\Theta]$, see \cite{Mi-ICM}.
Since the automorphy factor is affine linear the hypersurface $[\tilde\Theta]$
is periodic and thus descends to the quotient torus
$X$. We call the resulting hypersurface in $X$ the {\em theta divisor} $[\Theta]\subset X$.

\section{Abel-Jacobi Theorem}
In this section and throughout the rest of the paper we assume that $C$ is a {\em compact} tropical curve of genus $g$.
\subsection{Tropical Jacobian}
Denote by $\Omega(C)$ the space of global 1-forms on $C$. Given a collection of break points $(z_i,\eta_i)$ we can consider the values of forms on the $\eta_i$. This identifies $\Omega(C)$ with $g$ dimensional $\R$-vector space. The integral lattice $\Omega_\Z(C)\subset\Omega(C)$ consists of the forms taking integer values
on integer tangent vectors to the curve $C$.

Given a path $\gamma$ in $C$, any 1-form $\omega$ on $C$ pulls back to a piece-wise
constant classical (old-fashioned) 1-form on an interval. Thus one can define the integral $\int_\gamma \omega \in \R$.

Let $\Omega(C)^*$ be the vector space of $\R$-valued linear functionals on
$\Omega(C)$. Then the integral cycles $H_1(C,\Z)$ form a lattice $\Lambda$ in
$\Omega(C)^*$ by integrating over them. As in the classical complex geometry we
define the {\em Jacobian} of the curve $C$ to be
$$J(C):= \Omega(C)^*/H_1(C,\Z)\cong \Rg/\Lambda.$$
Note that the space $\Omega(C)^*$ is naturally endowed with the
(tautological) $\Z$-affine structure: the lattice $\Omega_\Z(C)^*\subset \Omega(C)^*$ is identified with the integer valued functionals on $\Omega_\Z(C)$. A set of break points gives a basis in $\Omega_\Z(C)$.

\begin{figure}[htbp]
   \includegraphics{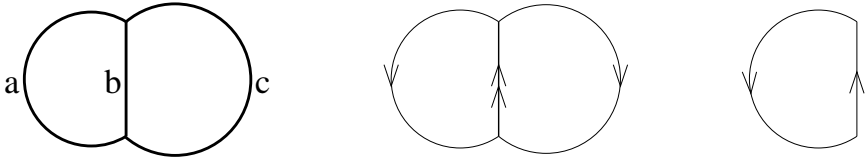}
   \caption{Two cycles $\gamma_1$ and $\gamma_2$ and their pairing
$Q(\gamma_1,\gamma_2)=a+2b$.}
   \label{fig:Q}
\end{figure}

The metric on $C$ defines a symmetric bilinear form $Q$ on the space of paths in
$C$ by setting $Q(\ell,\ell) := \length(\ell)$ for a simple (i.e., not
self-intersecting) path $\ell$ and extending it to any pair of paths
bilinearly, see Fig.~\ref{fig:Q}.
Since $C$ is a 1-dimensional simplicial complex any 1-cycle homologous to zero
is trivial at the simplicial chain level.
Therefore, $Q$ descends to a symmetric bilinear form on $H_1(C,\Z)$.

\begin{lemma}
The form $Q$ is positive-definite.
\end{lemma}
\begin{proof}
Since $C$ is a simplicial complex any 1-chain is given as a sum
$\sum_E a(E)E$, where $a(E)\in \Z$ and $E$ runs over all edges of $C$.
By bilinearity we have $$Q(\sum_E a(E)E,\sum_E a(E)E)=
\sum_E a(E)^2l(E)>0,$$
where $l(E)$ is the length of the edge $E$.
\end{proof}

Thus $Q$ provides an
isomorphism $\tilde Q: \Omega(C)^* \to \Omega(C)$. Under this isomorphism the
lattice $\Lambda$ is mapped isomorphically to the integral forms
$\Omega_\Z(C)\subset \Omega(C)$.  Thus the form $Q$ makes $J(C)$ into a
principally polarized tropical Abelian variety.

\begin{remark}
A version of the tropical Jacobian (as well as the tropical Picard group)
for finite graphs was introduced in \cite{Nagnibeda2} and \cite{Nagnibeda}
even before the explicit appearance of Tropical Geometry. From a tropical viewpoint
a finite graph is a graph defined over $\Z$ (i.e. with integer edge-lengths).
The finite graph Jacobian can be interpreted as the integer points
in the tropical Jacobian. The corresponding Abel-Jacobi Theorem was established in
\cite{Nagnibeda2}.
\end{remark}

\subsection{The Abel-Jacobi Theorem}
Once and for all let us fix a reference point $p_0\in C$. Given a divisor
$D=\sum a_i p_i$ we choose paths from $p_0$ to $p_i$. Integration along these
paths defines a linear functional on $\Omega(C)$:
$$\hat\mu(D)(\omega)=\sum a_i \int_{p_0}^{p_i} \omega. $$
For another choice of paths the value of $\hat\mu(D)$ will differ by an element
in $\Lambda$. Thus, we get a well-defined tropical analog of the {\em
Abel-Jacobi map} $\mu:\Div^d(C)\to J(C)$.
\begin{figure}[htbp]
   \includegraphics{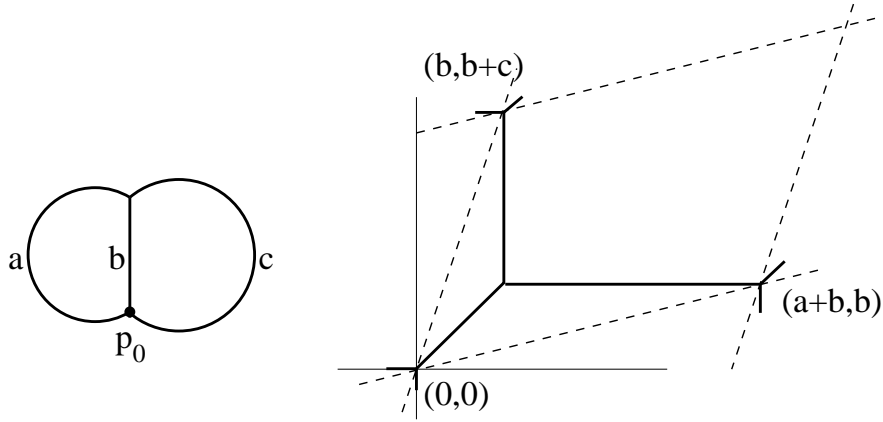}
   \caption{$\mu(C)$ in the tropical Jacobian $J(C)$.}
   \label{fig:jacobian}
\end{figure}
Note that the Abel-Jacobi map $\mu$ does not depend on the choice of a base
point $p_0$ if the degree $d$ is zero.
The dependence on $p_0$ for $d\neq 0$ will reappear later in solution to the Jacobi
inversion.

\begin{theorem}[Tropical Abel-Jacobi]\label{abel-jacobi}
For each degree $d$ the map $\mu$ factors through $\Picd(C)$:
\begin{center}
\mbox{}\xymatrix{ \Divd (C)\ar[dr]_{\mu} \ar[r] & \Picd (C) \ar[d]^{\phi}\\
     & J(C)}
\end{center}
so that $\phi$ is a bijection.
\end{theorem}
\begin{proof}
Here we only prove Abel's part: $\mu$ factors through and $\phi$ is injective.
Theorem \ref{jacobi_inversion} gives an explicit solution to the Jacobi
inversion.

Let $D=\sum p_i - \sum q_i$ be a divisor of a rational function $f$.
The gradient field $\nabla f$
 defines a linear functional on $\Omega(C)$ via integration along its integral
trajectories $t_j$ (counted with multiplicities). For some choice of paths from
$p_i$ to $q_i$ this functional coincides with $\hat\mu(D)$. On the other hand,
for any $\omega\in\Omega_\Z(C)$ we have
$$\sum_{j} \int_{t_i} \omega = \int_{\tilde Q^{-1}(\omega)} df=0.$$
Thus, $\mu(D)=0$ in $J(C)$.

Conversely, given any $D=\sum p_i - \sum q_i\in\Div^0(C)$ with
$\mu(D)=0$ we can choose $k$ paths $\ell_i$ on $C$, such that each
$\ell_i$ originates at $p_i$ and ends at $q_i$ and for any
$\omega\in\Omega(C)$ we have $\sum_{i} \int_{\ell_i} \omega=0$. Now
we can define a rational function $f(x)$ on $C$ by
choosing a path $\ell(x)$ from $p_0$ to $x$:
$$f(x):=\sum Q(\ell_i,\ell(x)).$$
Since $\sum Q(\ell_i,\gamma)=\sum \int_{\ell_i} \tilde Q(\gamma)=0$ for any
closed loop $\gamma$ in $C$ the function $f(x)$ is independent of the choice of
$\ell$. By construction $(f)=\sum p_i -\sum q_i$.
\end{proof}

\subsection{Jacobi inversion}
Given a tropical map $\phi:C\to X$ we can pull back any meromorphic
function on $X$ to a meromorphic function on $C$, see e.g.
\cite{Mi-ICM}. Let $L$ be a line bundle, $s$ be its
section and $D$ be the corresponding divisor.
Pulling back local representatives of $s$ defines a section $\phi^* s$ of some line bundle on $C$.
The divisor $\phi^* D$ of $\phi^* s$ depends only on $D$ and not on the choice  of $s$ (cf. \cite{Mi-ICM}
and \cite{AR}).

On the other hand, another section of $L$ defines a linearly equivalent divisor $D'\sim D$ on $X$, which in turn pulls back to a linearly equivalent divisor $\phi^* D'\sim \phi^* D$ on $C$. Thus the pull back line bundle $\phi^*L$ on $C$ is well defined. Alternatively, to define $\phi^*L$ one can pull back a defining \v{C}ech cocycle of $L$.

\begin{lemma}
The Abel-Jacobi map $\mu:C\to J(C)$ is tropical.
\end{lemma}
\begin{proof}
We need to show that for any point $p\in C$ there is a local chart $U\subset
R^{k-1}$ such that $\mu$ is the restriction of a $\Z$-affine map
$\R^{k-1}\to\R^g$. But this is clear because affine linear coordinates near a
$k$-valent vertex will also provide local coordinates in $\Omega(C)^*$ modulo
some {\em linear} relations among cycles.
\end{proof}

Next we will prove a refined version of the residue formula \eqref{eq:degree} for curves in
$\R^g$.

\begin{lemma}
Let $U$ be a connected open set of a curve tropically embedded in
$\R^N$ with boundary $\partial U =\{z_1,\dots,z_k\}$ and
$\nu_1,\dots,\nu_k$ -- the corresponding outward primitive tangent
vectors. Let $f$ be a rational function on $U$ with divisor $(f)=\sum
a_j p_j$. Then
\begin{equation}\label{eq:divisor}
\sum a_j p_j=\sum_{i=1}^k \left(\frac{\partial f}{\partial \nu_i}(z_i) \cdot
z_i-f(z_i) \cdot \nu_i \right),
\end{equation}
where the summation and equality takes place in the vector space $\R^N$, and
$p_j,z_i$ and $\nu_i$ are viewed as elements there.
\end{lemma}
\begin{proof}
The formula is additive with respect to gluing pieces and holds for affine
functions. On the other hand, at $z=p_i$ the statement is the definition of
$a_i$.
\end{proof}

For $\lambda\in \C^g$ let $\Theta_\lambda(x):=\Theta(x-\lambda)$ denote the translated theta function
and let $[\Theta_\lambda]$ be its divisor on $J(C)$ (the
corresponding line bundle $L_\lambda$ is different, but defines the same polarization as $L$).
 Let $D_\lambda:=\mu^*[\Theta_\lambda]$ denote the pull back of
$[\Theta_\lambda]$ to the curve via the Abel-Jacobi map $\mu:C\to J(C)$.
Note that both $[\Theta_\lambda]=[\Theta]+\lambda$ and $D_\lambda$ are well defined for $\lambda\in J(C)$.

\begin{theorem}[Jacobi Inversion]\label{jacobi_inversion}
For any $\lambda\in J(C)$ the divisor $D_\lambda$ is effective of degree $g$. There exists a universal $\kappa\in J(C)$ such that $\mu(D_\lambda)+\kappa=\lambda$ for all $\lambda\in J(C)$.
\end{theorem}
\begin{proof}
$D_\lambda$ is clearly effective since the theta function locally pulls back to regular sections. To calculate its degree we choose a set of break points $z_i, \eta_i$ in the interiors of edges and disjoint from the support of $D_\lambda$. Let $T\subset \hat C$ be the associated fundamental domain, and let $T^0$ denote its interior and $\bar T$ -- its closure.

Then $\partial T^0$ consists of  $g$ ordered pairs of break ends: $z_i^+\in T$ and $z_i^- \in \bar T\setminus T$. The vectors $\eta_i$ at  the $z_i$ provide an integral basis for $\Omega_\Z(C)^*$.

\begin{figure}[htbp]
   \includegraphics{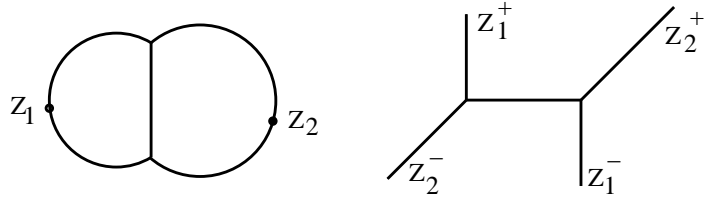}
   \caption{The interior of the fundamental domain $T$.}
   \label{fig:covering}
\end{figure}

Let $\hat\mu: \hat C\to \Rg$ denote the lifting of the Abel-Jacobi map $\mu$.
The path in $\bar T$ from $z_i^-$ to $z_i^+$ defines a cycle $\gamma_i$ in $C$,
that is $\hat\mu(z_i^+)=\hat\mu(z_i^-) + [\gamma_i]$ in $\R^g$. Differentiating the
quasi-periodicity equation of Lemma \ref{theta_period} gives
$$ d \Theta_\lambda(\hat\mu(z_i^+))- d \Theta_\lambda(\hat\mu(z_i^-))=\tilde Q
(\gamma_i)\in \Omega_\Z(C).$$
Note that the left hand side is independent of a lift of $\lambda$ to $\R^g$. The key observation is that $\tilde Q (\gamma_i)$ form a basis of $\Omega_\Z(C)$ dual to $\{\eta_i\}$, that is
$\<\tilde Q (\gamma_i), \eta_j\>=\delta_{ij}$. The residue formula  (\ref{eq:degree}) applied to the pullback $\hat\mu^*  \Theta_\lambda$ on $T^0$
gives $\deg D_\lambda=\sum_{i=1}^g 1= g$.

The second statement follows from the refined residue formula (\ref{eq:divisor})
applied to the theta function restricted to $\hat\mu(T^0)\subset \R^g$. Let $\nu_i=\hat \mu (\eta_i)$. Then we have
\begin{multline*}
\hat\mu(D_\lambda)= \sum_{i=1}^g
\left( \frac{\partial\Theta_\lambda}{\partial\nu_i}(\hat\mu(z_i^+)) \cdot
\hat\mu(z_i^+) - \Theta_\lambda(\hat\mu(z_i^+)) \nu_i \right) \\
-\left(\frac{\partial\Theta_\lambda}{\partial\nu_i}(\hat\mu(z_i^-)) \cdot
\hat\mu(z_i^-) - \Theta_\lambda(\hat\mu(z_i^-)) \nu_i \right).
\end{multline*}
Differentiating the right hand side with respect to $\lambda$ gives
$$\sum_{i=1}^g \left(  d \Theta_\lambda(\hat\mu(z_i^+)) \otimes \nu_i -  d
\Theta_\lambda(\hat\mu(z_i^-)) \otimes \nu_i\right) =
\sum_{i=1}^g \tilde Q (\gamma_i) \otimes \nu_i,$$
which is the identity element in $\mathrm{End}(\R^g)$. Thus, passing to the
quotient $J(C)$ we get $\mu(D_\lambda)=\lambda+\const$.
\end{proof}

Combining the Jacobi Inversion with the Abel-Jacobi Theorem we get the following statement.
\begin{corollary}\label{cor:jacobi}
Given a divisor $D$ of degree $d$ there is an effective divisor $D_\lambda=\mu^* [\Theta_{\mu(D)+\kappa}] $
of degree $g$ linearly equivalent to $D+(g-d)p_0$.
In particular, a degree $g$ divisor has a {\em canonical}
(independent of the base point $p_0$) effective representative in its class
of linear equivalence.
\end{corollary}

\subsection{Schottky problem and Torelli theorem}
The space of all principally polarized tropical Abelian varieties is the same as
the space of symmetric positive definite matrices modulo some discrete
automorphism group action. So its dimension is $\frac{g(g+1)}{2}$. On the other
hand the space of tropical curves of genus $g\geq 2$ is $3g-3$. Hence the
Jacobians form a subset $\mathcal{J}$ of positive codimension for $g\geq 3$
inside the space of all Abelian varieties. Description
of this subset is the tropical counterpart of the classical Schottky problem.

\begin{figure}[htbp]
   \includegraphics{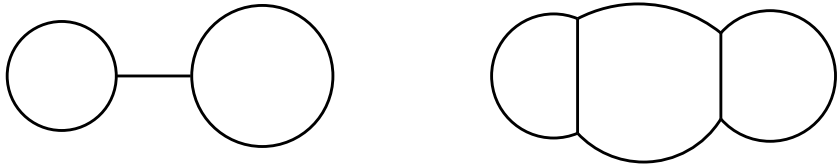}
   \caption{Counterexamples to Torelli in genus 2 and 3.}
   \label{fig:torelli}
\end{figure}

The na\"ive statement of the classical Torelli theorem fails for tropical
curves. For instance, in Fig. \ref{fig:torelli} the polarized Jacobian does not
see the length of the connecting edge in the genus 2 curve.
However we conjecture that following reformulation holds: the map from the
moduli space of tropical curves to $\mathcal{J}$ is tropical of degree 1.

\section{The Riemann-Roch Theorem}
In this section we continue to assume that $C$ is compact of genus $g$ and
look at the space of effective divisors realizing the same
class in the Picard group $\Pic(C)$.
Recall that for any divisor $D$ of degree $d$
we associate the linear system $|D|$ consisting
of all effective divisors linearly equivalent to $D$.
By Theorem \ref{abel-jacobi} $|D|$ is the inverse image of a point under
the Abel-Jacobi map $\Sym^d C \to J(C)$. Thus $|D|$ has the structure
of a compact (in topology induced from metric on $C$) CW-complex.

It is easy to see that (unlike the classical situations) the complex $|D|$ is not
of a pure dimension (cf. \cite{GK}, Example 1.11).
Nevertheless, we may speak of its dimension $\dim|D|$.
Namely, the following defintion was suggested by Baker and Norine.
\begin{definition} [\cite{BN}]\label{BNdef}
Let $D$ be an effective divisor on $C$.
The {\em dimension} $\dim |D|$ of the linear system $|D|$
is the maximal integer $r\geq 0$ with the following property:
for any effective divisor $R$ of degree $r$ the space $|D-R|$ is non-empty.
We set $\dim |D|=-1$ when $|D|=\emptyset$.
\end{definition}

The definition immediately implies that for effective $D$ and $D'$ one has
\begin{equation}\label{eq:sumdegree}
\dim|D|+\dim|D'| \leq \dim|D+D'| \leq \dim|D'|+\deg D.
\end{equation}
Plugging $D'=\emptyset$ into the second inequality gives us $\dim|D| \leq d$.
By the Jacobi inversion every divisor of degree $d\geq g$ is linearly equivalent to an effective one.
Again, directly by Definition \ref{BNdef} this implies that
\begin{equation}\label{eq:rre}
\dim|D|\geq d-g.
\end{equation}

If \eqref{eq:rre} turns into equality
then $D$ is called {\em regular}.
Otherwise we define the discrepancy number
$$\rho(D):= \dim|D| - d+g.
$$
An effective $D$ with $\rho(D) >0$ is called  {\em special}.

\begin{lemma}\label{lemma:2g-1}
Every divisor $D$ of degree $d\geq 2g-1$ is regular, i.e. $\dim|D| = d-g$.
\end{lemma}
\begin{proof}
Suppose $\dim|D|=s > d-g$.
Consider a degree $d-s$ divisor $D'$ with $|D'|\ne\emptyset$
(it exists since the Abel-Jacobi map is not surjective
on $\Sym^k C$ for  $k < g$).
It can be written as $D'=D-(D-D')=D-S$,
where $S$ is of degree $s \geq g$ and hence $|S|\ne\emptyset$.
Consequently, $|D'|\ne\emptyset$ which is a contradiction.
\end{proof}

The divisor
\begin{equation}\label{eq:canonical}
K:=\sum_{p\in C} (\val(p)-2)p,
\end{equation}
is called the {\em canonical divisor} of $C$. Here $\val(p)$ is the valence of $p$. By Euler's formula $\deg K=2g-2$. It is effective unless $C$ has one-valent vertices.

\begin{remark} In the degeneration picture a $k$-valent vertex $p$ of $C$
corresponds to the Riemann sphere with $k$ punctures. This explains
why the canonical divisor is defined as above:  the number of zeros
(the multiplicity of $p$ in $K$) of a rational form on  $\PP^1$ is
the number of poles (valence of $p$) minus 2.
In \cite{Mibook} all Chern classes are defined (in a similar combinatorial way)
for higher-dimensional tropical varieties.
\end{remark}

Since $C$ is a curve we may avoid dealing with higher-dimensional
cohomology by postulating the Serre duality. By introducing the cohomological
notation $h^0(D)= \dim |D|+1$ we set $h^1(D)=h^0(K-D)$.

\begin{theorem}[Tropical Riemann-Roch, cf. \cite{BN}, \cite{GK} and \cite{MZ1}]\label{thm:rr}
$$\dim|D|- \dim|K-D|= d-g+1.$$
\end{theorem}

The first version of our preprint \cite{MZ1} had a rather lengthy and involved independent
proof of this statement by means of tropical geometry.
However now we've learned of a much simpler proof of this statement
(in the case of finite graphs) by Baker and Norine \cite{BN} which
makes an elegant use of the chip-firing games technique (see e.g. \cite{BLS}, \cite{Biggs1}).
Since \cite{BN} appeared before \cite{MZ1} and contains a simpler proof
we see no point of pursuing our initial cumbersome argument here.
Below we just reproduce a geometrically adapted version of \cite{BN}.

There is already a couple of other papers that observed that the restriction
of \cite{BN} to the case of finite graphs is not essential. In particular,
\cite{GK} (where this observation appeared first) formally deduces
Theorem \ref{thm:rr} from \cite{BN}, and \cite{HKN}
adapts the argument of \cite{BN} to the case of metric graphs.
Note that one useful advantage of the approach of \cite{BN}
is the explicit description of all non-effective divisors
of degree $g-1$ (Corollary \ref{moderators}).

The main idea is to look for a canonical representative of a divisor
in its linear equivalence class.

Let us fix a point $p\in C$.
Any effective degree $d$ divisor on $C\setminus\{p\}$ may be considered just
as an unordered $d$-tuple of points. Thus to each such
divisor we may associate a non-decreasing sequence of distances
from these points to $p$. Given two divisors we may compare the
corresponding sequences in the lexicographic order (inserting some
number of zeroes in front of one of the sequences to make them
equal length).

\begin{definition}
We say that a divisor $D$ on $C$ is {\em $p$-reduced}
if its restriction to $C\setminus\{p\}$ is effective and
minimizes the distance to $p$ among such in the equivalence class of $D$.
\end{definition}

\begin{proposition}\label{p-unique}
For any fixed $p\in C$ there exists a unique $p$-reduced representative in every class in $\Pic(C)$.
\end{proposition}
\begin{proof}
Given $D$ of degree $d$ consider the minimal $m\in\Z$ such that $|D+mp|\ne \emptyset$.
Such $m$ is finite since it is obviously bounded from below by $-d$ and from above by $g-d$ by Jacobi Inversion.
Note that any effective $D'\sim D+kp$ for $k>m$ and with support disjoint form $p$ has strictly larger distance to $p$ than any element of $|D+mp|$ because of its higher degree.
Thus it suffices to minimize the distance to $p$ among elements of $|D+mp|$. Compactness of  $|D+mp|$ guarantees existence.

To show uniqueness let us suppose
that $D$ and $D'$ have the same distance sequence to $p$ and
there exists a {\em non-constant} rational function $f$ such that
$D'=D+(f)$. Since $C$ is compact the function $f:C\to\R$ has a maximum $M$.
Let $F_+=f^{-1}(M)$ be the locus of points where $f$ reaches its maximum.
Clearly $F_+$ is a subgraph of $C$ (possibly consisting of a single point).
Interchanging the r\^oles of $D$ and $D'$ if needed we may
assume that $F_+$ is disjoint from $p$.

All the boundary points of $F_+$ are the poles of $f$ of some order, thus these points are contained in $D$ with some positive multiplicities. Denote the divisor formed by them with $D_M$.
The points of $f^{-1}(M-\epsilon)$ can be enhanced with natural multiplicities
equal to the slope of $f$ on the corresponding edges of $C$. The resulting divisor
$D_{M-\epsilon}$ is linearly equivalent to $D_M$. In the same time, since $p\notin F_+$
the distance from $D_{M-\epsilon}$ to $p$ is smaller then the distance from $D_M$ to $p$.
Therefore $D$ cannot be $p$-reduced.
\end{proof}

\begin{remark}
The minimal $m$ argument for existence also shows that if $|D|\ne\emptyset$ then the $p$-reduced form is effective.
\end{remark}

\begin{definition}\label{defKplus}
A divisor $K_+$ on $C$ is called a {\em moderator} if there exists an acyclic orientation
on $C$ (i.e. a presentation
of $C$ as a 1-dimensional simplicial complex with a choice of an acyclic orientation
on the edges) such that
\begin{equation}\label{eq:K_+}
K_+=\sum_{p\in C} (\val_+(p)-1)p,
\end{equation}
where $\val_+$ stands for the number of outgoing edges.
\end{definition}

Analogously we can define $K_-:=\sum_{p\in C} (\val_-(p)-1)p$ by counting {\em incoming} edges,
which is also a moderator for the reversed orientation.
The notation is justified by the following proposition.
\begin{proposition}
The degree of $K_+$ and $K_-$ is $g-1$, and $K_+ + K_-=K$.
\end{proposition}
\begin{proof}
By definition we have
$$K_++K_-=\sum_{p\in C} (\val_+(p)-1 +\val_-(p)-1)p=\sum_{p\in C} (\val(p)-2)p=K.$$
To calculate the degree we observe every edge and every vertex (with negative sign) enters
in the coefficients of
either $K_+$ or  $K_-$ exactly once.
\end{proof}

\begin{lemma}\label{modnoteff}
$|K_+|=\emptyset$ for any moderator $K_+$.
\end{lemma}
\begin{proof}
Suppose that $K_+ + (f)$ is effective for some rational function $f$.
Take the same maximal locus subgraph $F_+\subset C$ as in the proof of Proposition \ref{p-unique}.
The presentation of $K_+$ as a moderator induces an acyclic orientation of the graph $F_+$.
Thus there must exist a sink vertex $q\in F_+$. If $f$ is locally constant in the neighborhood
of $q$ then it is also a sink in $C$ and $q$ enters $K_+ + (f)$ with a negative coefficient.
Otherwise $q$ enters in $(f)$ with the coefficient $-m$ where $m$ is not smaller than the number $n$
of edges of $C\setminus F_+$ adjacent to $q$. In the same time the coefficient of $q$ in $K_+$
is not greater than $n$.
\end{proof}

\begin{lemma}\label{prop:rr}
Given a divisor $D$ on $C$ exactly one of the following two holds.
Either $|D|\ne\emptyset$, or there is a moderator $K_+$ such that $|K_+-D|\ne\emptyset$.
\end{lemma}
\begin{proof}
By Proposition \ref{p-unique} we may assume that $D$ is $p$-reduced.
Consider a presentation of $C$ as a graph (with no loop edges) whose set of vertices $\{v_j\}$
contains $p$, the support points
of $D$ and points $C$ of valence greater than $2$. To choose an acyclic
orientation it suffices to order the vertices $v_j$. We do it inductively
starting with $v_0=p$.

Suppose that the first $k$ vertices are already chosen. Let us look at the edges $E_j^{(k)}$
connecting these vertices with the remaining vertices.
Note that $p$-reducibility implies that one of the remaining vertices, say $v$, enters $D$
with the coefficient smaller than the number of edges $E_j^{(k)}$ adjacent to it.
Otherwise we could move all the outer endpoints of all $E_j^{(k)}$ by some distance $\epsilon>0$
towards $p$ (this does not change the linear equivalence class)
and obtain a contradiction to $p$-reducibility. We set $v_{k+1}=v$.

Orienting each edge towards a smaller vertex gives a moderator $K_+$ with $K_+-D$ effective,
except possibly at $p$. If  $|D|=\emptyset$ then $D$ is not effective and it must have a negative coefficient at $p$ and thus we get $K_+-D$ effective at $p$.
If both $|D|$ and $|K_+-D|$ were non-empty we would get a contradiction with Lemma \ref{modnoteff}.
\end{proof}

This proposition has a notable corollary that also provides
the statement converse to Lemma \ref{modnoteff}.
\begin{corollary}\label{moderators}
Let $D$ be a divisor of degree $g-1$. If $|D|=\emptyset$ then $D$ is linearly equivalent to a moderator. Moreover, if in addition $D$ is $p$-reduced, then $D$ is a moderator.
\end{corollary}
\begin{proof}
By Lemma \ref{prop:rr}
there exists a moderator $K_+$ with $|K_+-D|\ne\emptyset$.
But the degree of $K_+-D$ is zero. Thus $D \sim K_+$.

If $D$ is $p$-reduced, then $K_+-D$ is effective, hence trivial.
\end{proof}

\begin{corollary}\label{cor:rr}
If $D$ is a divisor of degree $g-1$ then
$|D|=\emptyset$ if and only if $|K-D|=\emptyset$.
\end{corollary}
\begin{proof}
If $D=K_+$ is a moderator then $K-D=K_-$ is also a moderator.
\end{proof}

\begin{corollary}\label{div+1}
If $d<g-1$ and $|D+q|\neq\emptyset$ for every point $q\in C$
then $|D|\neq\emptyset$.
\end{corollary}
\begin{proof}
By Lemma \ref{prop:rr} if $|D|=\emptyset$ then there exists a moderator
$K_+$ with $|K_+-D|\ne\emptyset$. Since $\deg(K_+-D)>0$ there exists a point
$q$ such that $|K_+-D-q|$ is still non-empty. Thus $|D+q|=\emptyset$.
\end{proof}

\begin{proof}[Proof of Riemann-Roch theorem]
Because of Lemma \ref{lemma:2g-1} the Riemann-Roch holds for $d<0$ and $d> 2g-2$.
By symmetry between $D$ and  $K-D$ it suffices to prove only the inequality
\begin{equation}\label{RRE}
\dim|K-D| \geq \dim|D| - d+g-1
\end{equation}
for $0\leq d\leq 2g-2$. The inequality is trivial for $\rho(D)=0$.
Also applying (\ref{eq:rre}) to $K-D$ shows that (\ref{RRE}) holds if $|D|=\emptyset$.
So from now on we assume that $D$ is special.

We need to show that given any effective $R$ of degree $\rho(D)-1$ we have $|K-D-R|\neq\emptyset$.
Replacing $D$ with $D+R$ reduces the statement to showing that $|K-D|\neq\emptyset$ for any special
$D$ of degree $d\geq g-1$.

The case $d=g-1$ follows from Corollary \ref{cor:rr}.
If $D$ is a special divisor with $d\ge g$ we have $|D-R|\neq\emptyset$ for any effective $R$ of degree $d-g+1$,
and hence, $|K-D+R|\neq\emptyset$ by Corollary \ref{cor:rr}. Repetitive application of Corollary
\ref{div+1} shows that $|K-D|\neq\emptyset$.
\end{proof}

\section{Riemann's Theorem}
By taking the sum in the group $J(C)$ we may extend the Abel-Jacobi
map to $\mu:C\times\dots\times C\to J(C)$. This extension is still a tropical map.
Furthermore, since taking the sum is commutative
we also have a map $\mu:\Sym^k(C)\to J(C)$ (note though that $\Sym^k(C)$ is {\em not}
a tropical variety in the usual sense for any $k>1$).

The tropical counterpart of Riemann's theorem identifies
the subset $W_{g-1}=\mu(\Sym^{g-1}C)\subset J(C)$ with a translate of the theta divisor.
We establish it as a corollary of Jacobi's inversion.
We start with a couple of statements which themselves may be of independent interest.

\begin{definition}
The support $\supp|D|$ of the linear system $|D|$ is the set of points $q\in C$ such that $|D-q|\ne\emptyset$.
\end{definition}


Let $\Gamma$ be a subgraph of $C$. We denote by $D_{\partial\Gamma}$ its boundary divisor, that is the divisor consisting of the boundary point of $\Gamma$.

\begin{lemma}\label{lemma:support}
Let $\Gamma$ be a proper connected subgraph of $C$ and
$D_b$ be any break divisor of $\Gamma$ (cf. \S \ref{break_points}).
Then $\Gamma\subset \supp |D_b+D_{\partial\Gamma}|$.
\end{lemma}
\begin{proof}
Let $T_\Gamma$ be a fundamental domain of $\Gamma$ associated to some choice of break points constituting $D_b$
(as usual, we may identify $T_\Gamma$ with $\Gamma$).
By taking the closure of one of the connected components of $T_\Gamma \setminus D_{\partial\Gamma}$
instead of $\Gamma$ if needed we may assume that the multiplicity of any boundary point $v$ in
$D=D_b+D_{\partial\Gamma}$ is equal to the valence of $v$ in the graph $\Gamma$.
Indeed, $\val_\Gamma (v) >  \deg D_b \vert_v$,
since $v$ can support at most $\val_\Gamma (v) -1$ break points.
On the other hand, if $\val_\Gamma (v) \geq \deg D_b \vert_v+2$
then $v$ is not an endpoint of $T_\Gamma$ and it
breaks $T_\Gamma$ into several components.

Let $\epsilon>0$ be the distance between $\partial\Gamma$
and the vertices of $\Gamma\setminus\partial\Gamma$.
Here we assume presentation of the graph $\Gamma$ such that its set of vertices contains the support points of $D$ and points of $\Gamma$ of valence greater than $2$.
Consider the divisor $D'$ obtained from $D$ by moving the points $D\vert_{\partial\Gamma}$ inside $\Gamma$ by distance $\epsilon$, and define the proper subgraph $\Gamma'\subset\Gamma$ by by removing swept out edges: $\Gamma':=\{x\in \Gamma \suchthat \dist(x,\partial\Gamma) \geq \epsilon \}$.

Clearly, the boundary divisor of $\Gamma'$ is contained in the moved part of $D'$. On the other hand, the interior break points from $D_b$ in $\Gamma$ form a set of break points for $\Gamma'$. Hence $D'-D_{\partial\Gamma'}$ contains a break divisor on $\Gamma'$. The lemma follows now by induction replacing $\Gamma$ and $D$ with $\Gamma'$ and $D'$.
\end{proof}

Let $D$ be a degree $g$ divisor on $C$. We set $\lambda=\mu(D)\in J(C)$.
As in Corollary \ref{cor:jacobi} we denote with $D_\lambda=\mu^*[\Theta_{\lambda+\kappa}]$ the canonical effective form of $D$.

Recall that the theta divisor $[\Theta_{\lambda+\kappa}]$ defines the Voronoi cell decomposition of $\R^g$.
Let $Q_0$ be the interior of a maximal cell. We build a fundamental domain $Q\supset Q_0$ in $\R^g$
by including relative interiors of some boundary cells to satisfy the following ``boundary connectedness'' criterion:
if a cell $\tau$ is included then so is any cell of larger dimension incident to $\tau$.
\begin{figure}[htbp]
   \includegraphics{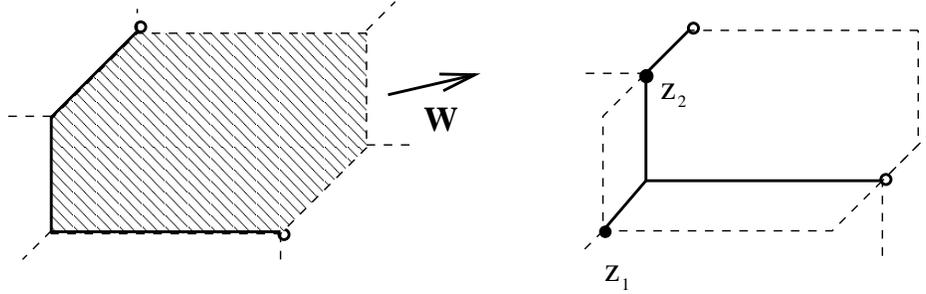}
   \caption{The fundamental domains $Q,T$ and the break divisor $D_\lambda=z_1+z_2$.}
   \label{fig:D_lambda}
\end{figure}
One way to choose such $Q$ is to take a generic linear functional $W$ on $\R^g$ and include those boundary cells whose baricenters are minimal with respect to $W$ in their corresponding equivalence classes (see Fig. \ref{fig:D_lambda}).

Furthermore, given a face $\sigma$ in $Q$ we can choose $W$ such that, in addition, the converse to the above criterion holds for $\sigma$ to the best possible extent. Namely, for any face $\tau\subset\sigma$ its class representative $\tau'$ in $Q$ is a face of $\sigma$. On Fig. \ref{fig:D_lambda} this holds for the leftmost face. In such case we say that $Q$ is {\em $\sigma$-complete}.

As before we denote by $\hat \mu:\hat C \to \R^g$ the lifting of the Abel-Jacobi map. Let $G=[\pi_1(C),\pi_1(C)]\subset \pi_1(C)$ be the commutator subgroup of the fundamental group of the curve, and let $\tilde C = \hat C/G$ be the intermediate abelian covering space.
By definition $\hat \mu$ factors thought the abelian covering $\tilde\mu: \tilde C \to \R^g$.

Let $\tilde T= \tilde\mu^{-1}(Q)\subset \tilde C$.
By the Abel-Jacobi theorem $\tilde T$ contains no cycles, hence it can be continuously lifted to a fundamental domain $T \subset \hat \mu^{-1}(Q)\subset \hat C$ of the curve.  Let $D_T$ be the associated pseudo-break divisor given by (\ref{eq:pseudo_break}).

%

\begin{lemma}\label{lemma:D_lambda}
$T$ is connected and $D_T$ is a break divisor. Furthermore, $D_T=D_\lambda$.
\end{lemma}
\begin{proof}
Suppose $z\in C$ enters in $D_T$ with multiplicity $k$. Consider the continuous lift of a neighborhood of $z$ to $\hat C$ such that $\hat z\in T$. Let $\eta_{i_1},\dots,\eta_{i_k}$ be the break vectors, that is the primitive tangent vectors at $\hat z$ which are {\em not} in $T$.
Also choose lifts of $\kappa$ and $\lambda$ to $\R^g$ such that $\Theta_{\lambda+\kappa}=0$ on $\bar Q_0$. By our boundary connectedness criterion for $Q$ the vectors $\eta_{i_1},\dots,\eta_{i_k}$ are precisely those primitive tangent vectors at $\hat z$ which are mapped outside of $\bar Q_0$. Then as we saw in the proof of the Jacobi inversion theorem $\Theta_{\lambda+\kappa}$ must have slope 1 into the directions $\hat \mu (\eta_{i_1}),\dots,\hat \mu (\eta_{i_k})$ and 0 along other primitive vectors.
Thus, the local pull back section $\mu^*\Theta_{\lambda+\kappa}$ has degree $k$ at $z$, that is $z$ enters in $D_\lambda$ with multiplicity $k$, and consequently, $D_T=D_\lambda$. Finally, Lemma \ref{lemma:pseudo_break} implies connectedness of $T$.
\end{proof}

Lemma \ref{lemma:D_lambda} identifies $\hat T:=\hat \mu^{-1}(Q)\subset \hat C$ as a disjoint union of fundamental trees on which $G$ acts faithfully and transitively by deck transformations. In fact, the components of $\hat T$ are fairly sparse in $\hat C$ in the following sense.
\begin{lemma}\label{lemma:covering}
A connected fundamental domain $T'$ in $\hat C$ intersects at most one component of $\hat T$.
\end{lemma}
\begin{proof}
Suppose $T'$ intersects two components $T_1,T_2 \subset \hat T$. Then there is a path $\hat \gamma$ in $T_1\cup T' \cup T_2$ from $q_1\in T_1$ to $q_2\in T_2$, with $q_1$ and $q_2$ projecting to the same point  in $\tilde C$. Let $\gamma$ be the image loop of $\hat\gamma$ in $C$. By considering a subpath of $\hat \gamma$ if needed we may assume that $\hat \gamma$ is simple and no two interior points of $\hat \gamma$ in $T_1 \cup T_2$ project to the same point in $\gamma$.

Consider a point $\hat q \in \hat \gamma$ which lies in $T'\cap(T_1 \cup T_2)$. Let $q\in C$ be its image in $\gamma$. Then $\gamma$ passes only once through $q$. This is impossible since $\gamma$ is trivial on the chain level in $C$.
\end{proof}


\begin{theorem}\label{supp|D|}
$q\in\supp|D|  \Longleftrightarrow \mu(q) \in [\Theta_{\lambda+\kappa}]$.
\end{theorem}

\begin{proof}
Let us choose the fundamental domains $Q$ and $T$ as above.

($\Longrightarrow$)
Let a point $q$ be in $C_0:=\mu^{-1}(J(C)\setminus [\Theta_{\lambda+\kappa}])\subset C$
(note that $C_0$ is not empty by our assumption)
and let $\hat q$ be its lift in $T\subset \hat C$. We can orient the infinite tree $\hat C$ such that $\hat q$ is the sink (the root) and every other vertex has exactly one outgoing edge. Thus every leaf of $T\subset \hat C$ is oriented inward. By identifying $T$ with $C$ we get an orientation on $C$, which by (\ref{eq:K_+}) defines a degree $g-1$ divisor  $K_+$. Each open leaf of $T$ contributes the corresponding break point into $K_+$. Thus, by Lemma \ref{lemma:D_lambda} we have $K_+=D_\lambda-q$.
Next we show that the described orientation on $C$ is acyclic. Then $D_\lambda-q$ is a moderator and, consequently, $q\not\in\supp|D_\lambda|=\supp|D|$ by Lemma \ref{modnoteff}.

Suppose there is an oriented cycle $\gamma$. Then its lift $\hat \gamma$ in $T\subset \hat C$ is a (connected) path. Let $v_1\in T$ and $v_2\in \bar T \setminus T$ be the two ends of $\hat \gamma$.
Then the path from $v_2$ to $\hat q$ in the tree $\hat C$ goes through $v_1$.
Note that $\hat \mu (v_1)$ lies in the relative interior of some boundary cell, say $\sigma_1$, of $Q$. Hence $\hat \mu (v_2)=\hat \mu(v_1)-[\gamma]$ lies in $\sigma_2=\sigma_1-[\gamma]$, a boundary cell of $\bar Q_0$ not included in $Q$. By continuity, there is a cell $\tau\succ\sigma_2$ which contains some part of $\hat\mu (\hat \gamma)$.

Then for a suitable choice of the linear functional $W'$ the new fundamental domain $Q'\subset \bar Q_0$ contains $\sigma_2$ and $\tau$ but does not contain $\sigma_1$. We choose the associated fundamental tree $T'\subset \hat \mu ^{-1} (Q')$ to contain $\hat q$. Let $v'_2\in T'$ be the lift of $\hat \mu (v_2)$, and let $\hat \gamma'\subset \hat T$ be the lift of $\gamma$ whose missing end is $v'_2$. Then $T'$ contains $\hat q$, $v'_2$ and some part of $\hat \gamma'$, but it does not contain $v'_1$, the other end of $\hat \gamma'$. Hence $\hat\gamma'$ is oriented opposite to $\hat \gamma$, and thus has to lie in a component of $\hat T$ different from $T$. But this contradicts Lemma \ref{lemma:covering}.

($\Longleftarrow$)
Let $\sigma\subset Q$ be the boundary face (closed in $Q$) whose interior contains the lift of $q$.
We may assume that our choice of $Q$ is $\sigma$-complete. Consider the subgraph $C_\sigma\subset C$ which is the projection of $\tilde \mu^{-1}(\sigma)\subset \tilde C$ under the covering map $\tilde C \to C$. Since $Q$ is $\sigma$-complete the subgraph $C_\sigma$ is closed.

Let $v$ be a boundary point of $C_\sigma$ and let $\xi_1,\dots,\xi_k$ be
the exterior primitive tangent vectors at $v$. Consider the continuous
lift of a neighborhood of $v$ to $\hat C$ such that $\hat v\in T$, then $\hat \mu (\hat v)\in Q$.
Neither of vectors $\hat \mu(\xi_i)$ can be inside the cone over $\sigma$ at $\hat\mu (\hat v)$.
Hence not all of them can be inside the cone over $Q$,
otherwise it would violate the tropical balancing condition at $\hat\mu (\hat v)$.

Hence, every boundary point of $C_\sigma$ supports at least one break point of $C$ whose tangent vector is exterior to $C_\sigma$. Consequently, $D_\lambda - D_{\partial C_\sigma}$ contains a break divisor for $C_\sigma$ (we may replace $C_\sigma$ by its connected component if needed). Thus we can apply Lemma \ref{lemma:support} to conclude that $C_\sigma\subset \supp|D_\lambda|$.
\end{proof}

\begin{remark}
In the case $q\in [\Theta_{\lambda+\kappa}]$ the corresponding orientation on $C$ always contains an oriented cycle, that is $D_\lambda-q$ is {\em not} a moderator. Indeed, the last paragraph of ($\Longleftarrow$) shows that every boundary point of $C_\sigma$ has an outward oriented exterior edge, but there is no sink in $C\setminus C_\sigma$.
\end{remark}

Note that if $\mu(C)$ intersect $[\Theta_{\lambda+\kappa}]$ in isolated points then $D$ is rigid, that is the linear system $|D|$ consists of a single element, namely $D_\lambda$. On the other hand, if $\mu(C)\subset [\Theta_{\lambda+\kappa}]$ then we get an analog of the classical result that the linear system $|D|$ is special. However, contrary to the classical situation, there are intermediate cases.

\begin{figure}[htbp]
   \includegraphics{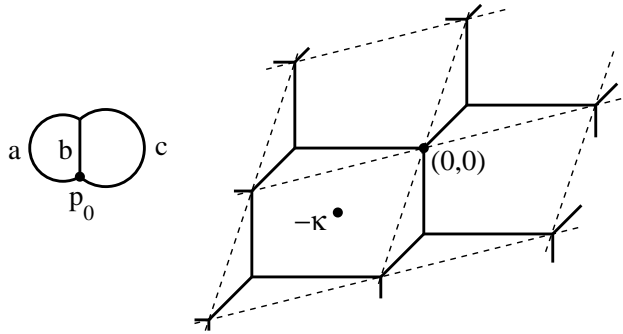}
   \caption{The theta divisor is the shift by $\kappa=(\frac{a+b}2,\frac{c+b}2)$
of the $W_1=\mu(C)$.}
   \label{fig:theta_divisor}
\end{figure}

\begin{corollary}[Riemann's Theorem]\label{thm:riemann}
$W_{g-1}+\kappa=[\Theta]$, where $\kappa\in J(C)$ is the Jacobi Inversion constant.
\end{corollary}

\begin{proof}
By Corollary \ref{cor:jacobi} a degree $g-1$ divisor $D$ is linearly equivalent to effective if and only if $|D_\lambda-p_0|\ne \emptyset$. By Theorem \ref{supp|D|} this is equivalent to $\mu(p_0)=0 \in   [\Theta_{\lambda+\kappa}] $. Hence we have
$$\lambda\in W_{g-1} \Longleftrightarrow
0 \in [\Theta_{\lambda+\kappa}] \Longleftrightarrow
0 \in [\Theta_{-\lambda-\kappa}] \Longleftrightarrow
\lambda \in [\Theta]-\kappa,
$$
where the second equivalence is because of $\Theta(-x)=\Theta(x)$.
\end{proof}

As a by-product we may now identify $\kappa\in\Pic^{g-1}$. By Corollary \ref{cor:rr} we have $W_{g-1}=\mu(K)-W_{g-1}$, and consequently
$$W_{g-1}+\kappa=[\Theta]=-[\Theta]=-W_{g-1}-\kappa=W_{g-1}-\mu(K)-\kappa.$$
But the theta line bundle defines a principal polarization, hence $[\Theta]$
cannot be stable under any non-trivial translation. Thus, $2\kappa=-\mu(K)$,
i.e. $\kappa$ can be interpreted as a tropical spin structure.

\bibliographystyle{alpha}
\bibliography{references}

\end{document}